\documentclass{cmslatex}
\usepackage[paperwidth=7in, paperheight=10in, margin=.875in]{geometry}
 \usepackage[backref,colorlinks,linkcolor=red,anchorcolor=green,citecolor=blue]{hyperref}
\usepackage{graphicx,graphics,caption, subcaption, enumerate,multirow, color}

\usepackage{amsmath,amssymb,amsfonts}
\usepackage{cases,bm,empheq}\usepackage{cite}
\usepackage{enumerate}
\renewcommand{\theequation}{\arabic{section}.\arabic{equation}}

\newcommand{\inpd}[2]{\left\langle #1, #2 \right\rangle}

\newcommand{\set}[1]{\left\{#1\right\}}

\newcommand{\eps}{\varepsilon}

\newcommand{\Real}{\mathbb {R}}

\DeclareMathOperator{\spn}{span}
\newcommand{\tr}{\textsf{T}}

\newcommand{\wt}[1]{\widetilde{#1}}
\newcommand{\Xb}{{\bar{X}}}

\newcommand{\C}{C}
\newcommand{\Cb}{\overline{C}}

\newcommand{\Lo}{\mathcal{L}}

\renewcommand{\d}{\ensuremath{\mathrm{d}}}
\newcommand{\dt}{ \ensuremath{\mathrm{d} t } }
\newcommand{\dx}{ \ensuremath{\mathrm{d} x} }
\newcommand{\dy}{ \ensuremath{\mathrm{d} y } }
 \newcommand{\Fd}{Fr\'{e}chet derivative}

\newcommand{\tm}{\wt{\rho}}

\newcommand{\FW} {Freidlin-Wentzell}

\newcommand{\e}{ \operatorname{E}}
 
 \newcommand{\var}{\operatorname{Var}}

  \newcommand{\hh}{\frac{1}{2}}
      \newcommand{\hfour}{\frac{1}{4}}

 \newcommand{\dxQ}{\nabla_x U}
 \newcommand{\dyQ}{\nabla_y U}
 
 \newcommand{\dxW}{\nabla_x W}

   \allowdisplaybreaks
\begin{document}
 \title{Multiscale  Gentlest Ascent Dynamics for Saddle Point
in Effective Dynamics of Slow-Fast System \thanks{Received date, and accepted date (The correct dates will be entered by the editor).}}

          \author{ 
 Shuting Gu\thanks{Department of Mathematics,
 City University of Hong Kong,Tat Chee Ave, Kowloon, Hong Kong SAR. shutinggu2-c@my.cityu.edu.hk.}
          \and  Xiang Zhou\thanks{Department of Mathematics,
 City University of Hong Kong,Tat Chee Ave, Kowloon, Hong Kong SAR. email: xiang.zhou@city.edu.hk.
The research of XZ  was supported by the grants from the Research Grants Council of the Hong Kong  SAR, China (Project No. CityU 11304314, 109113 and 11304715). ).}}

         \pagestyle{myheadings} \markboth{Multiscale   Gentlest Ascent Dynamics }
         {Shuting Gu and Xiang Zhou} \maketitle

          \begin{abstract}
             Here we present  a multiscale  method to calculate the saddle point
 associated with  the effective   dynamics arising from a stochastic  system which couples
 slow deterministic drift and   fast stochastic dynamics.
 This problem is motivated by the transition states on free energy surfaces in chemical physics.
  Our method is based on the    gentlest ascent dynamics
  which couples the position variable and the direction variable
  and has the  local convergence  to  saddle points.
  The dynamics of the direction vector is
  derived in terms of the covariance function with respective to the
  equilibrium distribution of the fast stochastic process.
 We apply  the multiscale numerical methods
to efficiently solve the obtained multiscale gentlest ascent dynamics,
{and  discuss the acceleration techniques based on the adaptive idea.}
The examples  of   stochastic ordinary and partial differential equations are presented.
          \end{abstract}


\section{Introduction }

The following slow-fast system  $(X(t), Y(t)) \in \mathcal{X}\times \mathcal{Y}$
 is a typical
dynamic system with  two disparate time scales:
\begin{subequations}{ \label{XY}}
 \begin{empheq}[left={ \empheqlbrace\,}]{align}
\dot {X}^\eps(t) & = f(X^\eps,Y^\eps) , \label{X}    \\
 \dot {Y}^\eps(t) &=  \frac{1}{\eps}b (X^\eps,Y^\eps)   + \frac{1}{\sqrt{\eps}} \sigma(X^\eps,Y^\eps) \eta(t),  \label{Y}
     \end{empheq}
 \end{subequations}
 where    $\eps$ is a   small parameter.
 $X^\eps$ is {the} slow variable and $Y^\eps$ is {the} fast variable.
For simplicity, we assume   $ \mathcal{X}=\Real^n$ and $\mathcal{Y}= \Real^m$. The functions  $f(x,y)=(f_1,f_2,\ldots, f_n): \Real^n \times \Real^m \to \Real^n$ and
$b(x,y)=(b_1,b_2,\ldots, b_m): \Real^n \times \Real^m \to \Real^m$ are two
smooth vector fields. The $m\times m$  diffusion matrix $\sigma(x,y)$  is assumed non-degenerate for all $x$ and $y$.
 $\eta(t)$ is the zero-mean Gaussian noise in $\Real^m$ with a certain covariance function $ \e(\eta(t)\eta(t'))$.
For example, $\e(\eta(t)\eta(t'))= \delta(t-t')$ means that $\eta$ is the white noise $\dot{W}$. We mainly consider this white noise case in this paper
for easy presentation, {though} the extension to the case with correlation structure is not difficult.
Note that in our model, the equation  \eqref{X} contains no diffusion term and the only random source comes from the $Y$-dependency of the function
$f$.

Many interests in   the above multiscale system concern the effective dynamics
of the slow variable, when  the
fast dynamics  can  be slaved and eliminated as the parameter $\eps$ tends to zero.  Indeed, for any fixed $T<\infty$, the slow variable in \eqref{XY},
 $X^\eps(t)$,   converges to a deterministic function $\Xb(t)$ satisfying the
 averaged equation
\[ \dot{\Xb}(t)=F(\Xb(t)),\]
during the time interval $t\in [0,T]$,  for  some function $F: \mathcal{X} \to \mathcal{X}$.
Under certain stronger conditions, this convergence is also uniform in $T$.
  This   is
the  typical result of  the averaging principle, which  has been
developed
in many classic mathematical literatures (see \cite{Khasminskii1966,Papanicolaou1977,KYSIAM2004}  and  \cite{PMM_book, FW2012}). On the application side,   one of the most important examples is the  {\it extended Lagrangian method} for  the
coarse-grained molecular dynamics (refer to \cite{ILP03,EVE2004} and \cite{PMM_book}).
In this application, the slow variable $X$ is usually the coarse-grained variables (collective variables) to describe the macroscopic features of the underlying complex system
 and the fast variable  $Y$
corresponds to the   microscopic variables of  full atomistic coordinates.
To map the full complex energy surface in $\mathcal{Y}$ into a low dimensional free energy surface in $\mathcal{X}$,
it is essential to solve  the averaged dynamics of the coarse-grained variable $\bar{X}$
with efficient numerical methods such as the heterogeneous multiscale method
(HMM) in \cite{EE2003,HMM07Rev,EVE2004}.
To faithfully recover the fluctuation for a finite $\eps$,
   \cite{KEVE2016}   used multiple replicas of fast process
in    HMM simulation.



%

 In this note we are interested in  the computation of  the saddle point in the effective dynamics $F$
 rather than the equilibrium fixed point.
 The exploration of the saddle point may help to
  investigate the phase space structure for $F$,  which has no closed form at all.
Specifically, we are interested in   index-1 saddle points
for the effective dynamics associated with the slow-fast system \eqref{XY}.
By index-1, we mean that the unstable manifold of the flow $F$
at the saddle point is one  dimensional.
The search of saddle point for the effective dynamics $F$ is not
an easy job in consideration of the lack of any analytical form of this dynamics.
Thus, any Newton-type method can not work in practice.
We shall use the gentlest ascent dynamics (GAD) developed in \cite{GAD2011}
to formulate  a new  dynamical system  associated with  $F$.
 This new system can be viewed as an {\it effective} dynamics for
a new multi-scale slow-fast system. To efficiently solve this new multi-scale system,
which is called MsGAD for convenience,  we apply the heterogeneous multiscale method or the seamless coupling method.

The original idea of the GAD is based on the
min-mode or  eigenvector-following methodology
(see, e.g., \cite{Crippen1971,Energylanscapes}).
The numerical developments based on  this methodology, such as
the dimer method(\cite{Dimer1999}), the activation-relaxation techniques (\cite{ART1998}), etc., have been
used for quite a few applications on    potential energy surface.  The applications of the GAD include   \cite{SamantaGAD,subcrit2010,GAD-DFT2015}.
The authors in
\cite{Bofill2014,Bofill2015} have  some further discussions on the properties of the GAD.
As an acceleration technique for the GAD,   a new
iterative formulation  and its algorithmic implementation for   gradient system  has
been  recently developed in \cite{IMF2014} and \cite{IMA2015}, respectively.
{The optimization-based idea is also implemented in \mbox{\cite{zhang2016sisc}}.}
One of the  advantages of the GAD is that its form  is  a continuous-time dynamic system, so it is very easy to fit in our current framework for the slow-fast dynamics
\eqref{XY}.  The application of GAD to sample the free energy surface has been
discussed  in \cite{samantaFES2014}     for temperature-accelerated
molecular dynamics.
{ The review of various saddle-point methods in computational material sciences
can found  in \mbox{\cite{npjRev2016}}.}


%
%
%
%
%

The saddle point problem is closely related to the transition-state calculation
for the randomly perturbed system with additive noise:  $ \dot{x} = F(x)+\sqrt{\eps}\dot{w}$.
The  optimal transition path in this additive noise case
is determined by  the quadratic form of the rate function in the traditional \FW\ large deviation theory(\cite{FW2012})
and many research works have shown that  the bottleneck on the optimal path is usually in the form of saddle points, especially
for the gradient system. So, many applications actually look for saddle points first rather than
compute the path directly.
 However,
the large deviation rate function for the slow-fast system \eqref{XY}
 has a much more complicated exponential form (\cite{Veretennikov200069,FW2012,Bouchet2016SlowFastLDP})
 and from a mathematical viewpoint,  little is known about the relationship  between the saddle point
and the transition state for the slow-fast dynamics.
In this paper, we do not intend to resolve the complicate issue of finding optimal  transition
 path for the slow-fast dynamics.
We  leave the efficient numerical scheme for the path calculation as the next project.
However,  from the viewpoint of  practical applications
of studying activated processes on free energy surface (\cite{ILP03,LaioParr2002,string-collective2006}),
  our {pursuit}  of the saddle point on
the effective dynamics $F$, which is the gradient of the free energy, is  of significance
if the transition state theory is used for rate calculation.
For example, the calculation of the nucleation rate from the bulk liquid to a crystalline solid
(e.g.~\cite{WRFJCP1996})
requires the location of  saddle point  on the free energy and the flux across the barrier at the saddle point.

 The rest of the paper is organized as follows.
 In Section \ref{ssec:ave} and Section \ref{ssec:GAD},  we review the averaging theorem for
 the slow-fast system \eqref{XY}
 and the gentlest ascent dynamics.
 In Section \ref{sec:MsGAD} we derive  the
 gentlest ascent dynamics for the averaged system, develop  the multiscale methods for computations
 and discuss  the algorithmic details, followed by
 a discussion on adaptive implementation and a remark on the connection to the central limit result.
  Section  \ref{sec:grad} is devoted to the gradient system
where  an extended potential energy function exists, in view of the practical importance
  of this class of models. Two examples are presented in Section \ref{sec:ex} to demonstrate
  our method and  the concluding remarks are given in Section \ref{sec:end}.

\subsection{Averaging principle of slow-fast stochastic dynamics}
\label{ssec:ave}
The averaging principle to derive the  effective dynamics of the slow-fast system \eqref{XY}
is based on  the ergodicity  assumption of the fast process $Y^\eps$.
Let $\wt{Y}^x$ be the solution of the equation
\begin{equation}\label{vY}
 \dot {\wt{Y}^x}  =  b (x,\wt{Y}^x)   + \sigma(x,\wt{Y}^x) \eta(t) ,
\end{equation}
 for any fixed parameter $x$. This is named as the virtually fast process.
 Assume that the virtually fast process $\wt{Y}^x$ is ergodic at every $x$ and
its unique invariant measure $\mu_x(\d y)$ has  a density function $\rho(x,y)$:
\begin{equation}
\mu_x(\d y ) = \rho(x,y)\d y=\frac{1}{Z(x)}   e^{- U(x,y)}   \d y,
\end{equation}
where the normalization factor $Z(x)$ is
\begin{equation} \label{def:Z}
Z(x) \doteq \int_{\mathcal{Y}} e^{-U(x,y)}\dy.
\end{equation}
\begin{remark}\label{rem:grad}
If  $b(x,y)= - \nabla_y U_1(x,y)$
and $\sigma(x,y)\equiv \sigma I$ for some potential energy function $U_1(x,y)$
and a constant $\sigma$, i.e.,
the fast dynamics is a gradient system,
then   $U(x,y) = \frac{2}{\sigma^2} U_1(x,y) $
and the equilibrium distribution $\rho(x,y) = e^{-\frac{2}{\sigma^2} U_1(x,y)}
/ \int  e^{-\frac{2}{\sigma^2} U_1(x,y)}\dy $.
In other cases,  we     simply set  $U(x)= -\log \rho(x,y)$ and $Z(x)\equiv1$.
\end{remark}

By the ergodicity assumption,    for any integrable function $u$, the expectation
with respect to $\mu_x$
can be estimated by the time average,
\[ \int u(y)\mu_x(\dy) = \lim_{T\to \infty} \frac{1}{T}\int_0^T  u(\wt{Y}(t))\dt.\]

The averaging principle  (cf. \cite{PMM_book,FW2012} and references therein)
states  that  as  $\eps\downarrow 0$,  the slow component of the system \eqref{XY}, $X^\eps(t)$, has a limit  $\Xb(t)$ satisfying the following ordinary differential  equation,
\begin{equation} \label{Xbar}
\dot{\Xb} = F(\Xb), ~~\mbox{ where }  F(x) \doteq \int f(x,y) \mu_x(dy).
\end{equation}

In most cases, the averaging function $F$ above has no closed formula,
and the solution   $\Xb$ has to be   approximated by
 the numerical methods.

\subsection{Gentlest Ascent Dynamics (GAD)}
\label{ssec:GAD}
For a dynamic system  $\dot{x}(t) = \varphi(x(t))$
where the flow   $\varphi$ is  $C^2$-smooth,
the gentlest ascent dynamics
locally converges to  saddle point of $\varphi$
by coupling the  position variable  and    the direction variable.
This dynamics, as a solution to the  saddle point problem,   can be viewed
as a counterpart   to the steepest descent dynamics for searching local minima.    In this note, we are only interested in the index-1 saddle point, i.e., the unstable manifold
of the saddle point is one dimensional.

The GAD  for the flow $\dot{x}(t) = \varphi(x(t))$ is the following extended system
for $(x,v,w)$,
\begin{subequations}\label{GAD}
\begin{empheq}[left=\empheqlbrace]{align}
  \dot{x}(t) &= \varphi(x) -2 \frac{\inpd{\varphi(x)}{ w }}{\inpd{w}{v} }v,\label{GAD-x}\\
  \gamma  \dot{v}(t) &=  D \varphi(x)v - \alpha v, \label{GAD-v}\\
 \gamma \dot{w}(t)&=  D\varphi(x)^\tr w - \beta w,  \label{GAD-w}
\end{empheq}
\end{subequations}
where  $\inpd{\cdot}{\cdot}$ is the inner product, $D\varphi(x)$ is the Jacobi matrix $\left(D\varphi\right)_{ij} \doteq
\frac{\partial \varphi_i}{\partial x_j}$.
Two scalars $\alpha$ and $\beta$ are Lagrangian multipliers to impose
certain  normalization conditions for $v$ and $w$. For instance,
if one uses the normalization condition  $\inpd{v}{v}=\inpd{w}{v}=1$, then   $\alpha=\inpd{v} {D\varphi(x)v}$ and $\beta=2\inpd{w}{D\varphi(x) v}-\alpha$.
By setting $\inpd{v(0)}{v(0)}=\inpd{v(0)}{w(0)}=1$,   the GAD flow
\eqref{GAD} for such choices of $\alpha$ and $\beta$ then will  preserve these two normalization equations. This technique to determine  $\alpha$ and $\beta$ will be applied  later
in many same  situations;
we  shall not  repeat the calculation of  these
Lagrangian multipliers.

The modified force in \eqref{GAD-x} has the effect of inverting the direction of
the original  force $\varphi(x)$ on the direction $v$ to sustain  ``ascent" flow  around the index-1 saddle point.
The dynamics
 \eqref{GAD-v} and \eqref{GAD-w}, if $x$ is frozen and  as time goes to infinity, tend to the
left-eigenvector $w$ and the right-eigenvector $v$
of the Jacobi $D\varphi(x)$ corresponding to the largest eigenvalue, respectively.
The coupling of $x$ and   $(v,w)$ is actually relaxed by
a finite positive  number  $\gamma$ in \eqref{GAD}, which
introduces a separation of   time scale artificially.
As $\gamma\downarrow 0$,   equation \eqref{GAD}
becomes  a two-scale   system where $x$ is   slow variable and $(v,w)$ are fast variables.
When the fast variables $v$ and $w$ have a single limit state as time goes to infinity, denoted by $v(x)$ and $w(x)$, respectively, $v(x)$ and $w(x)$ are the right and left eigenvector
of $D\varphi(x)$.
In this case,  at the limit $\gamma\downarrow0$,  the effective dynamics  of \eqref{GAD} is
 \begin{equation}
\label{L-GAD}
   \dot{x}(t) = \varphi(x) -2 \frac{\inpd{\varphi(x)}{ w(x) }}{\inpd{v(x)}{w(x)}} v(x).
   \end{equation}

The rigorous proof of the local convergence of the GAD \eqref{GAD}
to a nearby index-1 saddle point is presented in \cite{GAD2011} for any finite $\gamma$
and in the appendix of \cite{IMA2015} for the limit of vanishing  $\gamma$.

\begin{remark} \label{PDEGAD}
For the PDE case,  the transpose of the Jacobi matrix $D\varphi$  becomes
the adjoint of the variational derivative operator.  For instance, if the flow reads
$u_t = u_{xx} + u_x + f(u)$ with periodic boundary condition in 1-D, then the Jacobi matrix is $\partial_{xx}+\partial_x + f'(u)$
and its adjoint is $\partial_{xx}-\partial_x + f'(u)$.
\end{remark}

Since the Jacobi matrix $D\varphi$
is generally  asymmetric,
the  right direction $v$ and the left  direction $w$    are both  required in \eqref{GAD}
to   obliquely project the force $\varphi$ onto  $\spn\set{v}$ and $\spn\set{v}^\perp$,
except
 for  the gradient system  $\varphi(x)=-\nabla V(x)$, where the Hessian is symmetric.
The GAD for a gradient system involves only $v$:
\vspace{-10pt}
\begin{subequations}
\begin{center}
\begin{empheq}[left=\empheqlbrace]{align}
  \dot{x}(t) &= -\nabla V(x) + 2 \frac{\inpd{\nabla V(x)}{ v}}{\inpd{v}{v}} v,\label{GAD-g-x}\\
\gamma \dot{v}(t) & =   - \nabla^2 V(x)v + \inpd{v} {\nabla^2 V(x)v}v, \label{GAD-g-v}
 \end{empheq}
\end{center}
\label{GAD-g}
\end{subequations}
where $\nabla^2 V$ is the Hessian matrix of   potential energy function $V(x)$.
The   equation  \eqref{GAD-g-v}
is  in fact the steepest descent flow (rescaled by $\gamma$)
to minimize  the Rayleigh   quotient,   $ {\min}_{{\|v\|=1}} v^\tr H v$,
for the Hessian matrix $H\doteq \nabla^2 V (x)$.
So, the steady state $v(x)$ is actually the eigenvector of the Hessian $H$ at the lowest eigenvalue,
the so called ``min-mode".

In evolving the vector $v(t)$ in  \eqref{GAD-g-v},
it may not need the full Hessian matrix in practice.
 The computation of the Hessian-vector multiplication
$\nabla ^2 V(x) v$ is usually done by the finite difference method
as in the dimer method (\cite{Dimer1999,DuSIAM2012,zhang2014CiCP}),
because  this multiplication  is exactly the directional derivative along the $v$ direction:
$$\nabla^2 V(x) v = \frac{\d }{\d h} \nabla V(x+ hv)\vert_{h=0}
\approx h^{-1} (\nabla V(x+ hv)-\nabla V(x))
.$$

\section{ Multiscale   Gentlest Ascent Dynamics for Slow-Fast System}\label{sec:MsGAD}
\subsection{Formulation of GAD for slow-fast system}
We intend to extend the GAD to the slow-fast dynamics \eqref{XY}
to calculate    index-1 saddle points  of the effective flow $\dot{\Xb}=F(\Xb)$ defined in \eqref{Xbar}.
The  direct  application of the  GAD \eqref{GAD} to     equation
\eqref{Xbar}   gives
\vspace{-20pt}
\begin{subequations}
\begin{center}
\begin{empheq}[left=\empheqlbrace]{align}
  \dot{x}(t) &= F(x) -2 \frac{\inpd{F(x)}{ w }}{\inpd{v}{w}} v, \\
  \gamma  \dot{v}(t) &=  D F(x)v - \alpha v,  \\
 \gamma \dot{w}(t)&=  DF(x)^\tr w -\beta w,
\end{empheq}
\end{center}
\label{GAD-F}
\end{subequations}
where $F(x)=\int f(x,y)\rho(x,y)\dy$ by definition,
and  $DF(x)$ is the Jacobi matrix of
 $F(x)$. Note   the density $\rho(x,y)=Z(x)^{-1} e^{-U(x,y)}$.  We introduce
\begin{equation} \label{def:g}
g(x,y)\doteq - \nabla_x U(x,y),  \end{equation}
and
\begin{equation}\label{def:G}
G(x)\doteq \int g(x,y) \rho(x,y) \d y.
\end{equation}
By the definition of $Z(x)$ in \eqref{def:Z}, we simply see
\begin{equation} \label{def:G1}
\nabla_x  \log Z(x)  =Z^{-1}(x) \nabla_x Z(x)= \int g(x,y) \rho(x,y) \d y  =    G(x).
\end{equation}
We   calculate the Jacobi matrix as follows,
\[\begin{split}
(D F)_{ij}(x) & =\frac{\partial F_i}{\partial x_j} (x)
 = \frac{\partial }{\partial x_j} \left ( \int f_i(x,y) Z^{-1} (x) e^{- U(x,y)}   \d y \right )
\\
& = \int  \left ( \partial_{x_j} f_i(x,y) + f_i(x,y) g_j (x,y)  - f_i (x,y)
\partial_{x_j} Z(x)  Z^{-1}(x) \right)\rho(x,y) \d y
\\
& = \int  \left ( \partial_{x_j} f_i(x,y) + f_i(x,y) g_j (x,y)  - f_i (x,y)
G_j(x) \right)\rho(x,y) \d y
\\
& = \overline{ \partial_{x_j} f_i}(x) +\overline{ f_i g_j }(x)  -  F_i(x) G_j(x).
\end{split}
\]
To ease presentation, the overlined symbol $\overline{\theta}(x)$ for a bivariate  function $\theta(x,y)$ is
used  to define the
expectation with respective to $\mu_x(\dy)$, that is,
\[ \overline{\theta}(x)\doteq \int \theta(x,y) \mu_x(\dy).\]
So, $\overline{f}(x)=F(x)$ and $\overline{g}(x)=G(x)$ by this definition.

The Jacobi matrix of the effective dynamics
is thus given by \begin{equation}\label{DF}
DF (x)  = \overline{D_x f} (x) +  \Cb(x),
\end{equation}
where $D_x f(x,y)$ is the Jacobi matrix of $f(x,y)$ with respect to the variable $x$
and
\begin{equation}\label{def:C}
\C(x,y) \doteq  {f(x,y) \otimes g (x,y) }
-  F(x) \otimes G(x).
\end{equation}
The tensor product $u\otimes v$ for any two vectors $u$ and $v$ corresponds to  the matrix $[u_iv_j]$.

The  term  $\Cb(x)$ in \eqref{DF} comes from the
$x$-dependency of the equilibrium distribution $\mu_x(\dy)$.
$\Cb(x)$  actually is  the
covariance of $f$ and $g$ w.r.t. the distribution $\mu_x(\dy)$
because it is easy to verify that
$\Cb(x) = \overline{(f-F)\otimes (g-G)} (x)$.
 {This means that
one can use an alternative  form of \mbox{\eqref{def:C}}
as the follows}
\begin{equation}\label{def:C2}
{C(x,y) =  (f(x,y) -F(x))   \otimes ( g(x,y) - G(x)),}
\end{equation}
{if only the average quantity $\Cb$ is concerned.
One can also verify that the choice of  $C(x,y)=f\otimes g - F\otimes g$
also generates the same expectation $ {\Cb}$ as \eqref{def:C}.
}
\subsection{The Multiscale  GAD }
\label{ssec:MsGAD}

We shall address how to construct  multiscale schemes for the
system \eqref{GAD-F}.
We have obtained the expression of the Jacobi $DF$ in \eqref{DF},
which is   the ensemble average of the matrix
$D_x f + C$ w.r.t. $\mu_x$.
This important feature
allows us to view  the system \eqref{GAD-F}  as
an averaged equation of a   multiscale  system involving the original  fast variable $y^\eps$:
\vspace{-10pt}
  \begin{subequations}
\begin{center}
\begin{empheq}[left=\empheqlbrace]{align}
  \dot{x}^\eps(t) &= f(x^\eps,y^\eps) -2 \frac{\inpd{f(x^\eps,y^\eps)}{ w^\eps }}{\inpd{w^\eps}{v^\eps}} v^\eps,\label{GAD1-x}\\
   \dot {y}^\eps(t) &=  \frac{1}{\eps} b (x^\eps,y^\eps)   + \frac{\sigma(x^\eps,y^\eps)}{\sqrt{\eps}} \eta(t) ,\label{GAD1-y}\\
\gamma \dot{v}^\eps(t) &= \left(D_x f (x^\eps,y^\eps)  + \C(x^\eps,y^\eps) \right)v^\eps
  - \alpha^\eps v^\eps, \label{GAD1-v}\\
 \gamma\dot{w}^\eps(t)&=\left( D_x f (x^\eps,y^\eps) +   \C(x^\eps,y^\eps)  \right )^\tr  w^\eps -\beta^\eps w^\eps, \label{GAD1-w}
\end{empheq}
\end{center}
\label{GAD1}
\end{subequations}
where  the Lagrangian multipliers $\alpha^\eps$ and $\beta^\eps$ can be defined   as before
to enforce certain normalization conditons.
Here for any constant $\gamma$,
as $\eps\to 0$, the slow variables are $(x^\eps, v^\eps, w^\eps)$
and  the   fast variable is $y^\eps$.
We name the multiscale system like \eqref{GAD1} as the MsGAD (multiscale gentlest ascent dynamics).
Notice that the covariance term $\C$ is rank-$1$,
so the matrix-vector multiplication is actually
calculated simply as the
inner product,
for instance,  $\C(x,y)v = \inpd{g(x,y)}{v}f(x,y) - \inpd{G(x)}{v}F(x)$.

{
In the expressions of \eqref{def:C} or \eqref{def:C2},
{\it i.e.},
$\C(x,y)=f (x,y) \otimes g (x,y)
- F(x) \otimes G(x)$, or
$\C(x,y)= (f(x,y)-F(x)) \otimes (g(x,y)-G(x))$,
the averaged terms $F(x)=\bar{f}(x)$ and $G(x)=\bar{g}(x)$   already involve  the invariant measure $\mu_x$.
Hence,   equations \eqref{GAD-F} and  \eqref{GAD1} are not in the ``standard''    forms
like  the   equations \eqref{Xbar} and  \eqref{XY}
defined in Section \ref{ssec:ave}.
This issue will not bring any essential  challenges if the HMM is used
since     $F$ and $G$ can be estimated first from the sampling average (see below).
But in the  seamless coupling method to be introduced later in Section \ref{ssec:ss},
it is necessary to write the  system \eqref{GAD1}  in the  ``standard'' slow-fast
multiscale form as  equation \eqref{XY}.
It turns out that  a  single  fast variable $y$ is not sufficient:
one needs to  introduce another  process $\hat{Y}^x_t $ as an independent copy of the virtual fast process
$\wt{Y}^x_t$ to accomodate  the expectation   $F=\bar{f}$ (or $G=\bar{g}$).
Based on the fact that $\Cb=\overline{f\otimes g  - F\otimes  {g}}$,
we introduce
  \begin{equation}\label{def:C22}
  C_2(x,y,z) \doteq f(x,y)\otimes g(x,y) - f(x,z)\otimes g(x,y),
  \end{equation}
  then it is obvious that $\overline{C}(x) = \e_{y}\e_{z} \left[C_2(x,y,z)\right]$ where $\e_y$ and $\e_z$ are, respectively, the expectations w.r.t.  the independent random variables $y$ and $z$, which  both follow the same law $\mu_x$.
So, one   uses  the new    process $\hat{Y}^x_t $  to calculate the expectation   $F=\bar{f}$
w.r.t. $z$,
and the old $\wt{Y}^x_t$ to calculate
the  expectation  $\Cb=\overline{f\otimes g  - F\otimes  {g}}$ w.r.t. $y$  as
before.  This approach is based on the equivalent form of $\Cb$ :
\begin{equation*}\label{C4SS}
 \Cb(x)=  \int  \left[ f(x,y)\otimes g(x,y) -
 \left(\int f(x,z) \hat{\mu}_x(\d z) \right)\otimes g(x,y)    \right] \mu_x(\dy),
 \end{equation*}
where $\hat{\mu}_x(=\mu_x)$ is the equilibrium distribution of the iid copy $\hat{Y}^x$.

 However, for  the alternative  form in equation \eqref{def:C2}: $\C(x,y) =  (f(x,y) -F(x))   \otimes ( g(x,y) - G(x))$,
 $\Cb=\overline{(f-\bar{f})\otimes (g-\bar{g})}$ contains three expectations.
 One can run  three fast processes $y,z,w$ simultaneously as iid copies,
and then  letting $C_3(x,y,z,w)\doteq (f(x,y)-f(x,z)) \otimes (g(x,y)-g(x,w))$,
one has $\e_y\e_z\e_w \C_3 = \overline{\C}$.
It is interesting to find that it is possible to 
run   only  two processes $y$ and $z$ in this case by choosing 
$\hat{C}_2 (x,y,z)\doteq (f(x,y)-f(x,z)) \otimes (g(x,y)-g(x,z))$.
Then the calculation  shows
$\e_y \e_z \left[\hat{C}_2 (x,y,z)\right]=2(\overline{f\otimes g}-F\otimes G)=2\overline{\C}$.
Thus, one should use $\frac12 \hat{C}_2$ rather than $\hat{C}_2$
to estimate the covariance matrix $\Cb$.
In summary, there may be different   unbiased estimators of 
$\Cb$ if multiple streams of the fast processes are used.
The variances of these estimator could  be further analytically analyzed or numerically compared.
We do not further pursue this issue and simply  use the scheme based on the expression \eqref{def:C22} for numerical examples in this paper.
}

\subsubsection{HMM}
\label{sssec:HMM}
Now we discuss the framework of  the HMM ({\it heterogeneous multiscale method},
  \cite{EE2003,HMM07Rev}) for our averaged GAD system \eqref{GAD-F} and
  the MsGAD \eqref{GAD1}.
  There are two parameters $\eps$ and $\gamma$ in the
 multiscale GAD.
 We   let  $\eps$ tend to zero  in  \eqref{GAD1} to obtain the equation \eqref{GAD-F}
 by the   averaging principle.  We can also  further select
 a small $\gamma$ in \eqref{GAD-F} to obtain an equation like
 \eqref{L-GAD}.

The procedures of the HMM for the MsGAD are as follows.
Select a macroscopic  time step size $\Delta t$ for  evolving $x$
and $\Delta \tau$ for evolving $v$ and $w$
(usually $\Delta \tau=\Delta t$);
  select  a microscopic  time step size $\delta t$ for  evolving $y$.
The HMM scheme with forward Euler consists of the following steps.
\begin{enumerate}
\item Use the macro-solver
\[ x_{n+1} = x_n + \left(F_n - 2  v_n  \inpd{F_n}{v_n} / \inpd{w_n}{v_n} \right) \Delta t, \]
where  $x_n$ is the approximation value of $\Xb(n \Delta t)$
and $F_n$, $v_n$, $w_n$ are estimated below.
\item Apply the micro-solver with time step size $\delta t $  to $M$ micro-steps for the fast dynamics
\[  y_{n,m+1} = y_{n,m} + \frac{\delta t}{\eps} b(x_n, y_{n,m}) + \frac{\sigma(x_n,y_{n,m})}{\sqrt{\eps}}\sqrt{\delta t}\, \eta_{n,m},
\]
$m=0,1,2\ldots,M-1$. Here $\set{\eta_{n,m}}$ are iid $\mathcal{N}(0,1)$ random variables.
{The initial is choose as the  ``warm start'' :  $y_{n,0}\doteq y_{n-1,M}$.
Note that only the ratio $\delta t/\eps$ is needed here;
equivalently one can view this ratio as the true step size for the virtual
fast process $\tilde{Y}^x$ defined  in \mbox{\eqref{vY}}. }
\item Estimate $F_n$, $G_n$ and the Jacobi matrix $(DF)_n$:
\[\begin{split}
F_n &=\frac{1}{M} \sum_{m=1}^M f(x_{n}, y_{n,m}), ~~
G_n =\frac{1}{M} \sum_{m=1}^M g(x_{n}, y_{n,m}), \\
(DF)_n  &   =  \frac{1}{M} \sum_{m=1}^M \bigg (  D_xf ( x_n, y_{n,m})
+ f(x_n,y_{n,m})\otimes g(x_n,y_{n,m})\bigg )
\\  &~~~\quad - F_n \otimes G_n.
\end{split}
\]
\item Solve the right and left direction $v_n$ and $w_n$
for $K$ steps
by using the time step size $\Delta \tau$ and
the initial $v_{n,0}=v_{n-1}$, $w_{n,0}=w_{n-1}$:
\[
\begin{split}
\hat{v}_{n,k+1} &= v_{n,k}  +  \Delta \tau (DF)_n v_{n,k} ,
\quad \hat{w}_{n,k+1} = w_{n,k}  +  \Delta \tau (DF)^\tr_n w_{n,k} ,
\\
v_{n,k+1} &= \frac{ \hat{v}_{n,k+1} } { \lvert \hat{ v}_{n,k+1} \rvert},
\quad
w_{n,k+1} =\frac { \hat{w}_{n,k+1} }  {  \inpd{v_{n,k+1}}{\hat{w}_{n,k+1}}  },
\end{split}
\]
$k=0,1,2\ldots,K-1$. Then $v_n=v_{n,K}$ and $w_n=w_{n,K}$.
\end{enumerate}

The microscopic time step size has to be much smaller than $\eps$;
  $\delta t/\eps$ is the effective step size in solving the virtually fast process
$\wt{Y}$.
The time   length, $M\times\delta t$, should be sufficient long for the fast process
$Y^\eps$ to relax toward the
equilibrium distribution. The number of steps $M$  also should  be   large enough to reduce the
statistical error in estimating the averaged quantities   $F$, $G$ and $DF$.
The number of steps $K$ is to take care of a possible small  constant $\gamma$.
 $K=1$ actually works in principle for many cases. A larger $K$ can give
a better accuracy for eigenvector, {though} it also has more computational burden.
There is no requirement for $\Delta \tau$ as long as the ODE solver is stable;
 the time step size  $\Delta \tau$ can be simply set as $\Delta t$.

 \begin{remark}
On the choice of the macro-solver for $x$ and $v,w$,  an explicit scheme
with larger stability region is preferred, such as the stablized Runge-Kutta methods.
On the micro-solver for the virtually fast process, a numerical scheme  for SDE with higher order weak convergence rate for the {\emph {long}} time integration is preferred
to capture the equilibrium distribution better.
For instance, when $\sigma$ is a constant, one can use the stochastic
Heun method (\cite{kloeden2013}) which requires two force evaluations at each time step
or a non-Markovian scheme proposed recently in \cite{Leimkuhler2015}:
\[y_{n,m+1} = y_{n,m} + \frac{\delta t}{\eps} b(x_n, y_{n,m}) + \frac{\sigma\sqrt{\delta t}}{2\sqrt{\eps}}\, (\eta_{n,m}+\eta_{n,m+1}), ~~m=0,1,\ldots,\]
where $\set{\eta_{n,m}}$ are iid $\mathcal{N}(0,1)$ random variables.
 \end{remark}

\subsubsection{Seamless coupling method  (SCM)}
\label{ssec:ss}
The seamless coupling strategy proposed in \cite{Seam2009} does not need the back and forth communication
of the {macro}- and micro-states of the system.
It was thought as  a boosting algorithm: by increasing
the small parameter $\eps$ to a larger number $\eps'$ in  the   slow-fast system,
the seamless scheme simultaneously solves this boosted system  with the time step size $\wt{\Delta} t$ which is smaller than the macro-time step size $ \Delta t$ in the HMM.
It is potentially more efficient than the HMM if the micro-model is difficult.
Following this idea,
 we    increase $\eps$    in \eqref{GAD1}
to  a  larger number, say $\eps' =\eps\lambda$ for a constant $\lambda>1$.
As stated in previous section, to use \eqref{def:C2}, we need to introduce an independent
copy $z^\eps$ for the fast process to handle the double expectation in the covariance matrix $\Cb$.
In summary, the seamless coupling method solves   the following system:
\vspace{-10pt}
 \begin{subequations}
\begin{center}
\begin{empheq}[left=\empheqlbrace]{align}
  \dot{x}^\eps(t) &= f(x^\eps,y^\eps) -2 \frac{\inpd{f(x^\eps,y^\eps)}{ w^\eps }}{\inpd{w^\eps}{v^\eps}} v^\eps,\label{GAD2-x}\\
   \dot {y}^\eps(t) &=  \frac{1}{\eps\lambda} b (x^\eps,y^\eps)   + \frac{1}{\sqrt{\eps\lambda}} \eta_1(t) ,\label{GAD2-y}\\
     \dot {z}^\eps(t) &=  \frac{1}{\eps\lambda} b (x^\eps,z^\eps)   + \frac{1}{\sqrt{\eps\lambda}} \eta_2(t) ,\label{GAD2-z}\\
{\gamma} \dot{v}^\eps(t) &= \left(D_x f (x^\eps,y^\eps)  \right) v^\eps +
 f(x^\eps,y^\eps) \inpd{g(x^\eps,y^\eps)} {v^\eps}  \notag
 \\
 &~~\quad
 -f(x^\eps,z^\eps) \inpd{g(x^\eps,y^\eps)} {v^\eps}
  - \alpha^\eps v^\eps, \label{GAD2-v}\\
 {\gamma} \dot{w}^\eps(t)  &= \left(D_x f (x^\eps,y^\eps)  \right)^\tr w^\eps +
 g(x^\eps,y^\eps) \inpd{f(x^\eps,y^\eps)} {w^\eps}   \notag
 \\
 &~~\quad
 -g(x^\eps,z^\eps) \inpd{f(x^\eps,y^\eps)} {w^\eps}
  - \beta^\eps v^\eps, \label{GAD2-w}
\end{empheq}
\end{center}
\label{GAD2}
\end{subequations}
where $\eta_1$ and $\eta_2$ are two iid copies of $\eta$.
The joint fast processes $y^\eps$ and $z^\eps$ correspond to the
  equilibrium distribution $\mu_x(\dy)\times \mu_x(\d z)$.
  It is clear that as $\eps\to 0$, the effective dynamics
  of \eqref{GAD2} is just equation \eqref{GAD-F}.
The system \eqref{GAD2} is solved by any standard  ODE/SDE solver,
{such as the Euler scheme},
with a common time step size $\wt{\Delta }t$ for all five components.
  \cite{PMM_book} suggests  a  time step size $\wt{\Delta }t =\Delta  t /M$, where $M$ is the micro-step number previously used in the HMM.
After a sufficient long time for relaxation,
the solution $x^\eps(t)$ will go close to the index-1 saddle point
of the flow $F$, then one may  average $x^\eps(t)$ within a time interval
to improve the accuracy.

 {
 To end this subsection, we discuss the difference between the HMM and the seamless coupling
 method (SCM). One can refer to the book
\mbox{\cite{PMM_book}} and references therein for    general discussions.
For the application of these two methods to the search of stationary points here,
we emphasize some key and special differences below.
Firstly, the HMM and SCM solve two related but different
 systems:  equation \eqref{GAD-F}  and equation \eqref{GAD2}, respectively.
 The HMM solves the effective dynamics \eqref{GAD-F} by
simulating its  multiscale
version \eqref{GAD1} at the vanishing $\eps$ limit,
while the SCM  solves the pre-asymptotic system \eqref{GAD2}
with an effective parameter $\eps'=\eps\lambda>0$. So,
as $\eps'$ tends to $0$, the solution of the SC
tends to the solution of the HMM.
In other words, the result  of the SCM is an approximate solution
to the true saddle point of the effective force $F$
and the error is proportional to $\sqrt{\eps'}$ (see Section \ref{ssec:clt}).
In the long time run,
this also means that, with a fixed $\eps'>0$,
 the trajectory of the SCM
will eventually hop between multiple saddle points (if exist)
as the time goes to infinity.
Secondly,
in the HMM, the two expectations in  the covariance matrix $\overline{C}$
are computed via the sample average from a   single stream
of random samples from $\mu_x$; while in the  SCM,
this double expectation is realized  by two streams of
random samples for the virtual fast process.
Therefore, the  understanding of the seamless coupling here is
more than the classic {boosting} idea
since the classic seamless coupling for slow-fast two-scale system
does not have  to run two  independent streams of the  virtual fast process.
Thirdly,    the  HMM needs calculation of the sample average, which
has high computational costs but controls variance very well.
The SCM does not take sample averages at each iteration,
which is fast but is more noisy.
 }

\subsection{Adaptive strategy}\label{Sec2_3}

The quantity of our interest is the saddle point, not the whole
trajectory.
This focus on the final destination, regardless of the accuracy
of the trajectories at the early stage, facilitates the application  of
     adaptivity ideas to fine tune some parameters in our algorithms.
In order to reduce the variance of the solution
and improve the practical efficiency,
a very simple idea is to apply the SCM first  
and then switch to the HMM after an appropriate time period.
We generally discuss a few more useful methodologies below.

\begin{enumerate}
\item Average the output:   $\tilde{x}(t) \doteq \frac{1}{t-t_0}\int_{t_0}^{t} x(s)\d s $
after a period of    ``burn-in''  time    $t_0$.
This is the widely used technique in Monte Carlo simulation and the stochastic gradient dynamics
to reduce the variance of the original noisy output.
This principle works for both HMM (particularly for a small $M$) and the SCM without
any   extra computational burden, but for the SCM,
the effect of improvement  is more significant.

\item Choose the sample size $M$ adaptively in the HMM:
Use a small $M$ at the first stage for a quick search of the area
of the solution, then use a large $M$ to suppress the variance when the fluctuation starts to dominate the errors.
One might simply link the increasing rate of $M$ to the time, say $M\propto t^{r}$ for some $r>0$.

\item  Reduce the small parameter $\eps'$ adaptively in the SCM:
The errors of the seamless coupling method
depend on the small parameter of $\eps':=\eps\lambda$ and the step size $\wt{\Delta} t$.
A smaller $\eps'$ means a smaller deviation  of $x^\eps$ from the true solution
(associated with $F$).
So, one can  decrease the value of $\eps'$ as time runs.
Since the effective step size for the virtual fast processes
is $\wt{\Delta} t/\eps' $, then one also has to decrease
the step size  $\wt{\Delta} t$ so that $\wt{\Delta} t/\eps'$
 either is fixed or decreases too.
For example,  the typical decay of $\wt{\Delta} t$ at time $t$
in stochastic optimization
corresponds to the scaling   $ \wt{\Delta} t \propto t^{-1}$;
equivalently, the step size at $k$-th iteration
is  $\wt{\Delta} t_k =\wt{\Delta} t_1 /\sqrt {k}$.
Then the decay of   $\eps'$     may be  set as
$\eps'(t)\propto  t^{-1+c}$ for a  small constant $0\leq c< 1$,
so that $\wt{\Delta} t/\eps' \propto t^{-c}$.

\end{enumerate}

We shall demonstrate the effectiveness of the above three ideas
in the numerical section.

%

\subsection{Numerical calculation of the right direction  $v$ }
\label{ssec:re}

In the algorithms presented above,
we formally  use the full matrix of the Jacobi $DF$
in    the dynamics
for the right and left directions, $v(t)$ and $w(t)$, respectively.
If the dynamics is in the form of PDEs,
then as pointed out in Remark \ref{PDEGAD},  the Jacobi $DF$
and its transpose are just the variational derivative and its adjoint,
both of which are not difficult to derive in most cases.
However,  in some real applications such as atomistic models with long range interaction,
it is not efficient or even feasible to store this Jacobi matrix
element by element.
In what follows, we discuss this challenge
and   the numerical remedies   according to
the specific features of   slow-fast system.
 We start from the right direction $v$.
It is noted that in many cases,  the dimension of $\mathcal{X}$
is   much lower than the dimension of $\mathcal{Y}$. This means that
$D_x f(x,y)$ is a small scale matrix and its all entry-wise values can be
calculated and stored with reasonable cost.

 If the fast dynamics is of gradient type, i.e., $b(x,y) = -\nabla_y U_1(x,y)$ and $\sigma(x,y)\equiv \sigma I$,
then  $U=\frac{2}{\sigma^2}U_1$ and $g(x,y)=-\frac{2}{\sigma^2} \nabla_x U_1(x,y) $.
The covariance matrix  is $\C(x,y)=-\frac{2}{\sigma^2}  f(x,y) \otimes \nabla_x U_1(x,y)$.
Assume  the analytic form of $\nabla_x U_1(x,y)$ is available in practice, then
 the rank-1 matrix $C$ is easy to calculate:
 $C(x,y) v = -\frac{2}{\sigma^2}  f(x,y) \inpd{ \nabla_x U_1(x,y)}{v}$.

However, for non-gradient dynamics,     the  equilibrium distribution $\rho(x,y)$
 has no closed formula  available.
 For the dynamics of $v$ in \eqref{GAD1-v},
 the finite difference scheme
 \begin{equation}\label{DFv}
 \begin{split}
 (D F(x)) v & =
  \lim_{h\to 0} (F(x+hv) - F(x))/h \\
 &=\lim_{h\to 0} h^{-1} \int  f(x+hv,y)\rho(x+hv,y) - f(x,y)\rho(x,y)\dy \\
\end{split}
 \end{equation}
  can be applied to compute the matrix-vector multiplication $(D F(x)) v $.
If one knows  the closed form of the Jacobi matrix $D_x f$,
$(D_x f(x,y)) v $ thus can be directly computed by the matrix-vector multiplication.
The remaining term $(DF - D_x f)v = \Cb(x)v$
can be computed by the finite difference scheme.
 By  Remark \ref{rem:grad}
 and equations \eqref{def:g} and \eqref{def:G},
  we have  $g=\nabla_x \log \rho$ and
 $ \C(x,y) =f(x,y)\otimes \left( \nabla_x \log \rho(x,y) \right)
 =\rho^{-1}f(x,y)\otimes \left( \nabla_x  \rho(x,y) \right).$
Then,
 \[
  \begin{split}
   {\C}(x,y)v &= \rho^{-1}f(x,y)\inpd{v}{\nabla_x   \rho(x,y)}
   \\ &=
 \rho^{-1} f(x,y) \lim_{h\to 0} (  \rho(x+hv,y)-  \rho(x,y) ) / h.
  \end{split}\]
It follows that
\begin{equation}
\begin{split}
\Cb(x) v &= \int \C(x,y) v  \rho(x,y)\dy
\\
&=  \lim_{h\to 0} h^{-1}  \int f(x,y)(  \rho(x+hv,y)-  \rho(x,y) )   \dy.
\end{split}
\end{equation}
So,   \eqref{DFv} can be rewritten as
\begin{equation} \label{DFv2}
 (DF(x))v = \int D_x f(x,y) v  \rho(x,y) \dy + \lim_{h\to 0} \frac{1}{h}  \int f(x,y)(  \rho(x+hv,y)-  \rho(x,y) )   \dy.
 \end{equation}

In practice, a finite step size $h$ is used in \eqref{DFv} or \eqref{DFv2}.
Then, in the MsGAD,
besides the process $\wt{Y}^{x}_t$ slaved at $x$,
a second independent virtually fast process $\wt{Y}^{x+hv}_t$  is also required
 to obtain the  distribution $\rho(x+hv,\cdot)$.
  The second integration in \eqref{DFv2} is   computed from the time  averages
  from $\rho(x+hv,\cdot)$ and $\rho(x,\cdot)$.

\subsection{Numerical calculation of the left direction  $w$}

In many cases, we have to evolve  the dynamics for   $w$ direction in \eqref{GAD1-w},
except that the original slow-fast system \eqref{XY} is of the gradient type jointly (see
Section \ref{sec:grad} below).
It  is worthwhile to note   that  even if both the slow and  fast components
are of (parametrized) gradient types separately, but  the coupled  system  \eqref{XY}  is not jointly driven by a single potential energy function,  then the resulted averaging equation
\eqref{Xbar} might  still not be a gradient system and the Jacobi matrix
$DF(x)$ is {not}  symmetric.  Refer to Section \ref{eg:2D} for such an example.

If    the fast dynamics is of the gradient type,
as we mentioned above for $\C(x,y)v$ calculation,  then  it is also quite simple to calculate
  $\C(x,y)^\tr w$ by using $\nabla_x U_1(x,y)$.

Yet, when the fast dynamics is of the non-gradient type,  the matrix-vector multiplication
 $ ( DF(x))^\tr w $, in particular $C(x,y)^\tr w$,  will impose a very  severe computational challenge
 in this rather general situation, because it is not possible to apply the trick
 of  directional derivative.
We do not have any perfect solution for this non-gradient case
and the only method to evaluate the derivatives is the numerical evaluation.
For  $DF= \overline{D_x f + \C}$,
 we assume that  the full matrix  $D_x f(x,y)$ has an expression
  to compute and  its transpose
 $(D_x f(x,y))^\tr$ is obtained by a numerical transpose operation
 (or simply the adjoint operator in the infinite dimensional setting of the PDE case).
  Then the    term
  \[ \overline{\C^\tr  {w}} =  \int  \theta (x,y) \nabla_x   \rho(x,y) \dy,
  ~~~\mbox{ where }~~~   \theta(x,y) \doteq  \inpd{f(x,y)}{w} , \]
involves the sensitivity analysis for the equilibrium distribution $\nabla_x   \rho(x,y)$.
Although this sensitivity could be analytically derived in some case,
it can be approximated by the numerical derivative of $\rho$
  \begin{equation} \label{dxrho}
  \partial_{x_j} \rho(x,y) =  \lim_{h\to 0} (\rho(x+he_j,y)- \rho(x,y))/ h,
  ~~\quad j=1,\ldots,m,
  \end{equation}
  where $e_j$ is the unit vector along $x_j$ axis.
  Then the $j$-th component is
   \[ (\overline{\C^\tr  {w}})_j =   \lim_{h\to 0}  \int  \theta (x,y) (\rho(x+he_j,y)- \rho(x,y))/ h\, \dy.
\]
By choosing a finite number  $h\ll 1$,  this brute-force calculation  will have to add $m$ independent fast processes
  $\wt{Y}^{x+he_j}$, $j=1,\ldots,m$,
  to account for $\rho(x+he_j,\cdot)$.
The scheme \eqref{dxrho} is computationally feasible only when the dimension of $\mathcal{X}$
is low.  In the next subsection, we show that the $x$-derivative  of $\nabla_x\rho$
can be transformed to the $y$-derivative of some function by the perturbation
analysis for the equilibrium density $\rho(x,y)$. However, the numerical challenge
in this new form might still exist.

 \begin{subsection}{Connection with the central limit theorem}
 \label{ssec:clt}
 In the last part of this section,
  we shall give a remark on the connection with the normal deviation from the averaged system
 { and how this normal derivation theory can be useful in theory
  to estimate the fluctuation in the SCM. }
 In   the MsGAD, the  Jacobi $DF(x)$ is important: it determines
 the linearization of the averaged dynamics $F$.
 This linearization  also plays an important  role in
 approximating  $X^\eps-\Xb$, the difference of the solutions to
 the multiscale system
  \eqref{XY} and the averaged system  \eqref{Xbar}.
  Define the normalized difference
 \begin{equation}
 \label{1082}
 \xi^\eps_t \doteq \frac{1}{\sqrt{\eps}} (X^\eps_t-\Xb_t).
 \end{equation}
By Theorem 3.1, $\S 7.3$ in \cite{FW2012},   under the assumption of strong mixing,
as $\eps \to 0$,
the normalized difference converges weakly to the solution of the  following SDE
\begin{equation}\label{xi}
\dot{\xi}(t) = (DF(\Xb_t) )\, \xi (t)+ \beta(\Xb_t) \eta(t),
\end{equation}
where  $DF(\Xb_t)$ is the Jacobi matrix $DF$ at the averaged solution
$\Xb$ and the diffusion term
$\beta$ is the $m\times m$ matrix such that  $$\beta(x)\beta(x)^\tr =A(x)\doteq
\lim_{T\to \infty}\frac{1}{T}\int_0^T \int_0^T
  (f(x,\wt{Y}^x_s) - F(x))\otimes (f(x,\wt{Y}^x_t) - F(x)) \d t \d s ,$$
 which can be   estimated by running a
 long trajectory of the fast process $\wt{Y}^x$.

 If  the formal asymptotic expansion is applied to derive the
equation \eqref{xi},
the drift term in \eqref{xi} has a different form, denoted as
$B$,
\begin{equation} \label{B}
B(x)= \overline{D_x  f}(x) + \int  \left(  \int_0^\infty \nabla_y \e^y f(x, \wt{Y}^x_\tau)\d \tau
\right) D_x b(x,y)  ~\mu_x(\dy),
\end{equation}
where $\e^y$  is the expectation for the distribution of $\wt{Y}^x$ with initial
$\wt{Y}^x_0=y$.
Refer to the appendix in \cite{Bouchet2016SlowFastLDP} for this formula.
Comparing \eqref{B} with \eqref{DF}, we find that
$\Cb(x)$ should  be equal to  the double integral term on the right-hand side of \eqref{B}.
The proof of this fact is attached in the appendix of this paper.
Numerically, the formula \eqref{B} is not friendly due to the  differential operator $\nabla_y$
for the expectation term $\e^y$.

In summary, the normalized difference  $\xi$ in \eqref{xi} in the central limit theorem
shares the same drift flow as our dynamics for the right direction $v$ in \eqref{GAD-F}.
{Hence, the numerical methods we developed here for the MsGAD may be
useful to the calculation of $\xi$.}
The more interesting observation is that the above result
can aid in understanding the fluctuations in the seamless coupling method, at least in theory.

\begin{remark}
The above convergence theorem for $\xi^\eps_t \to \xi(t)$
as $t$ tends to infinity  
can also be applied to the effective GAD system \eqref{GAD-F}
and its multiscale system \eqref{GAD2}.
That is to define $\bar{X}$ and $X^\eps$    in \eqref{1082}
as the trajectories of the GAD system \eqref{GAD-F} (computed from the HMM)
and that of \eqref{GAD2} (computed from the SCM), respectively.
Then the corresponding normalized deviation satisfies the equation 
in form of  \eqref{xi} (after redefining the $F$ and $\beta$ terms accordingly
for the GAD rather than for the original dynamics),
which  characterizes the fluctuations
of the SCM.
One can furthermore linearize the equation \eqref{xi} around the   saddle point $x^*$
 to obtain a rough estimate of
the fluctuation around the true solution $x^*$.

\end{remark}

\end{subsection}
%

%
%
%

\section{The MsGAD for Gradient System}
\label{sec:grad}
As mentioned earlier,  one important example of the slow-fast system \eqref{XY}
 in practice
is the extended Lagrangian method for the coarse-grained molecular dynamics.
In this example, an energy potential $U(x,y)$ exists in the extended space  $\mathcal{X}\times\mathcal{Y}$ to drive the slow-fast system.
 For this gradient system, we shall see that the averaged system $\dot{\Xb}=F(\Xb)$
is also a gradient system, i.e., a (free energy) function $W(x)$ {exists} in $\mathcal{X}$ such that $F(x)=\nabla W(x)$.
The general discussions for the GAD above
can be {greatly simplified} much for this gradient case.
Now we assume the multiscale system \eqref{XY}   has the special form
\vspace{-20pt}
\begin{subequations}
\begin{center}
\begin{empheq}[left=\empheqlbrace]{align}
  \dot{X^\eps} &= -\dxQ(X^\eps,Y^\eps), \\
  \dot{Y^\eps} &= -\frac{1}{\eps} \dyQ(X^\eps,Y^\eps)  + \frac{1}{\sqrt{\eps}}\sigma\eta(t),
\end{empheq}
\end{center}
\end{subequations}
where
$\sigma$ is assumed a scalar constant and $\eta$ is a standard white noise.
$U(x,y)$ is a potential energy function for $(X^\eps,Y^\eps)$.
The equilibrium measure  of the virtually fast process  is
\begin{equation}
\mu_x(\d y) = \rho(x,y)\d y = \frac{1}{Z(x)} e^{-\frac{2}{\sigma^2} U(x,y)} \d y,
~~\ Z(x) \doteq \int e^{-\frac{2}{\sigma^2} U(x,y)} \d y.
\end{equation}
As $\eps \downarrow 0$, the averaged equation  for $\Xb$ is
\begin{equation} \label{Xbar1}
\dot{\Xb} = F(\Xb),
\end{equation}
where
$$
F(x) =  -\int \dxQ(x,y) \rho(x,y) \dy.
$$
By the definition of $Z(x)$, we can get
\[\begin{split}
\nabla_x  \log Z(x)
                    & = \frac{2}{\sigma^2} \int -\dxQ(x,y) \rho(x,y) \dy
                   = \frac{2}{\sigma^2} F(x).
\end{split}
\]
Thus, the averaged equation (\ref{Xbar1}) can be rewritten as a gradient system
\begin{equation}\label{eqn:aver}
\dot{\Xb} = -\dxW(\Xb),
\end{equation}
where the effective potential is
\begin{equation}\label{W} W(x) = -\frac{\sigma^2}{2} \log Z(x)=-\frac{\sigma^2}{2}
\log \left( \int e^{-2U(x,y)/\sigma^2} \d y \right).
\end{equation}
By the calculation of \eqref{DF},
the Hessian matrix of $W(x)$  is
\begin{equation}
\label{HessW}
\nabla_x^2 W(x) = - D F(x)
= \overline{\nabla_x^2 U}(x)  -  \frac{2}{\sigma^2} \overline{\dxQ \otimes \dxQ } (x) + \frac{2}{\sigma^2}
\overline{\dxQ}(x) \otimes \overline{\dxQ} (x).
\end{equation}
The right hand side contains   the Fisher information matrix
of the
invariant measure $\mu_x (=\rho(x,y)\dy) $:
$$-\e_{\mu_{x}}\left[\nabla_x^2 \log \rho\right ] = \frac{4}{\sigma^2}\left(
\overline{\dxQ \otimes \dxQ } (x) -
\overline{\dxQ}(x) \otimes \overline{\dxQ} (x)
\right).$$

\subsection{Example}
\label{ssec:grad-ex}
The following extended potential is
widely used in free energy sampling and the {coarse grained} molecular dynamics simulation
\cite{EVE2004}:
\begin{equation} \label{def:Q}
U(x,y) = V(y) + \hh \kappa\lvert x-q(y)\rvert^2,
\end{equation}
 where   $q=(q_1(y),\ldots, q_n(y))$ is a given function
 mapping fast variables in $\mathcal{Y}$ to the space  $\mathcal{X}$, which
 defines  { coarse-grained variables}. $\kappa > 0$ is a parameter coupling
 the potential of the microscopic system and the coarse-grained variables. Ideally, $\kappa$   should be infinitely large, but  it is   a large constant in practice.
The slow-fast dynamic system  associated with \eqref{def:Q} is
\vspace{-20pt}
\begin{subequations}
\begin{center}
\begin{empheq}[left=\empheqlbrace]{align}
  \dot{X^\eps} &= -\kappa(X^\eps-q(Y^\eps)),
  \label{ext-x} \\
  \dot{ Y^\eps} &= -\frac{1}{\eps} \Big (\nabla V(Y^\eps) -  \kappa( Dq(Y^\eps))^\tr (X^\eps-q(Y^\eps))\Big )+ \frac{1}{\sqrt{\eps}}\sigma \eta(t).
\end{empheq}
\end{center}
\label{extLag}
\end{subequations}
where $Dq(y)$ is the Jacobi matrix  $(\partial_{y_j}q_i)$.
 In order to sample the space $\mathcal{X}$ with the correct marginal
 equilibrium distribution  in the extended Lagrangian method.,
the slow dynamics \eqref{ext-x} actually should  also be independently driven by an
Brownian motion.
Since  here we   only study  the saddle point
rather than the distribution of $X^\eps$,
we   only concern  the deterministic steepest descent drift flow    in \eqref{ext-x}.

For this   example, we have that
$\nabla_x U(x,y)=\kappa (x-q(y))$ and
$\nabla_x^2 U(x,y)\equiv\kappa I$, where $I$ is the identity matrix.
It follows that
 $ \overline{\dxQ}=\kappa (\bar{x} -  Q(x))$
 and $\overline{\nabla_x^2 U}\equiv\kappa I$,
where
\[ \bar{x}=\int x \rho(x,y)\dy, ~~\quad Q(x) \doteq \bar{q}(x)= \int q(y) \rho(x,y)\dy, \]
and
\[\rho(x,y) =Z(x)^{-1}\, e^{-\frac{2}{\sigma^2}  V(y)}\, e^{-\frac{\kappa}{\sigma^2}  (x-q(y) )^2}. \]
Then,
the effective potential $W$ is
\[ W(x) = -\frac{\sigma^2}{2}
\log Z(x) =-\frac{\sigma^2}{2}
\log  \left( \int e^{-\frac{2}{\sigma^2}  V(y)} e^{-\frac{\kappa}{\sigma^2}  (x-q(y) )^2} \dy \right).
\]
and the effective force for the slow variable is
\[   F(x)   =-\nabla_x W(x) = -\kappa(\bar{x} -Q(x)).
\]
The Hessian matrix of $W$ for this example, by the result in \eqref{HessW}, is
\begin{equation}
\nabla_x^2 W(x)=\kappa I  -\frac{2\kappa^2}{\sigma^2}
\Big (
\overline{(x-q)\otimes (x-q)} (x)- (\bar{x}- Q(x)) \otimes (\bar{x}-  Q(x))
\Big),
\end{equation}
and hence
\[ (\nabla_x^2 W ) v =\kappa v  -\frac{2\kappa^2}{\sigma^2}\left(  \overline{\,\inpd{x-q}{v} (x- q)\,}- \inpd{\bar{x}-Q}{v} (\bar{x}-Q)\right).
\]
 The dynamics  for  direction    $v$ in the GAD,
$\gamma \dot{v}   =  - (\nabla^2 W)v+ \alpha v$, is
\[ \gamma \dot{v}   =
\frac{2\kappa^2}{\sigma^2}
\left( \overline{\,\inpd{x-q}{v}  (x-q)\,}- \inpd{\bar{x}-Q(x)}{v}(\bar{x}-Q(x)) \right) +  \alpha' v
\]
by absorbing the term $\kappa$ into the new Lagrangian multiplier $\alpha'$.

In summary,  for this example, the GAD for the averaged system is
\vspace{-20pt}
\begin{subequations}
\begin{center}
\begin{empheq}[left=\empheqlbrace]{align}
  \dot{x}(t) &= -\kappa(\bar{x}-Q(x)) +2 \frac{\inpd{\kappa(\bar{x}-Q(x)) }{ v }}{\inpd{v }{v}} v,
  \label{GADF-x}\\
  \gamma \dot{v}(t)  & =
\frac{2\kappa^2}{\sigma^2}
\left( \overline{\,\inpd{x-q}{v}  (x-q)\,}- \inpd{\bar{x}-Q(x)}{v}(\bar{x}-Q(x)) \right) +  \alpha' v.
\label{GADF-v}
\end{empheq}
\end{center}
\label{GADF}
\end{subequations}
If  $\kappa\to\infty$, then $\rho(x,y)\to \wt{Z}(x)^{-1} \delta(x-q(y)) e^{-\frac{2}{\sigma^2}V(y)}$, where $\wt{Z}(x)\doteq \int \delta(x-q(y)) e^{-\frac{2}{\sigma^2}V(y)}\dy$. And it follows that $W(x)\to \mathcal{F}(x)\doteq -\frac{\sigma^2}{2}\log \wt{Z}(x)$ (up to a constant),
where $\mathcal{F}$ is the free energy surface of the coarse-grained variable $x$.
For further details on the GAD for this type of question and the applications
in coarse-grained molecular dynamics simulation, refer to the publication \cite{samantaFES2014}.

\subsection{The dimer method  for gradient system}
The GAD for the gradient system \eqref{eqn:aver} requires the calculation of
the multiplication of   Hessian $\nabla_x^2 W$ and   vector $v$, which can be
numerically approximated by the finite difference scheme when the Hessian itself is difficult to obtain.
This idea has been widely used in the dimer method (\cite{Dimer1999,DuSIAM2012}).
The one-side  finite difference scheme is
\[ (\nabla_x^2 W )\,v \approx \frac{\nabla_x W(x+hv)-\nabla_x W(x)}{h}=
-\frac{F(x+hv)-F(x)}{ h}. \]
To evaluate $F(x+hv)$, which is $\int f(x+hv,y) \rho(x+hv,y)\dy$,  one simple  approach is to
run a second independent fast process $\wt{Y}^{x+hv}_t$ to obtain the
  density $\rho(x+hv,\cdot)$ parametrized at $x+hv$.
 When the central finite difference scheme is used, the third process for
 $x-hv$ is required. This approach  leads to an undesired extra  burden of
 simulating multiple fast processes.
 In contrast,
 by using the formula derived in  \eqref{HessW},
   only   one  trajectory for the virtually fast process is required.

\section{Numerical Examples}
\label{sec:ex}
To illustrate how the   MsGAD works, we present two numerical
examples below. The first  is a system consisting of
 a two dimensional  ordinary differential equation and
 a two dimensional stochastic ordinary differential equation.
It is not a gradient system. We apply the HMM
 (Section \ref{sssec:HMM}) in the MsGAD.
The second is a system of an Allen-Cahn partial differential equation and
a stochastic  Allen-Cahn partial differential equation.
This second system has an extended
potential functional. We apply and compare  both the HMM and the seamless coupling scheme (Section \ref{ssec:ss})
to this second example.

\subsection{A two-dimensional  example}
\label{eg:2D}
We consider the following system on $\mathcal{X}\times \mathcal{Y}=\Real^2\times\Real^2$
\vspace{-20pt}
\begin{subequations}
\begin{center}
\begin{empheq}[left=\empheqlbrace]{align}
  \dot{X}_i   &= - \sum_{j}D_{ij}X_j + Y_i^2,  \label{eg1-x}\\
 \dot{Y}_i & =   - \frac{1}{\eps}\frac{ Y_i}{ \Gamma_i(X)}  + \frac{1}{\sqrt{\eps}}\sigma \eta(t).\label{eg1-y}
 \end{empheq}
\end{center}
\label{eg1}
\end{subequations}
The vector field  $\Gamma(x)=(\Gamma_1(x), \Gamma_2(x)): \mathcal{X}\to \mathcal{Y}$
is  given.
$D=(D_{ij})$ is a $2\times 2$ symmetric matrix. $\sigma$ is a constant.

 The processes
$\set{Y_i(t)}$ are   independent  Ornstein-Uhlenbeck  processes  parametrized by     $X=x$.
The equilibrium distribution of $Y=(Y_1,Y_2)$  is
the product measure of   $ \mathcal{N}(0,  {\sigma^2\Gamma_i(x)}/{2})$. The calculation shows that  the limit equation has a  closed form
\begin{equation}\label{eg1-F}
\dot{\Xb}_i = - \sum_{j}D_{ij}\Xb_j + \frac{\sigma^2}{2}\Gamma_i(\Xb).
\end{equation}

Note that {if $D$ is positive-definite,} then \eqref{eg1-x} can be rewritten as
$\dot{X}=-\nabla_X ( \frac12 X^\tr D X - \sum_i{Y_i^2 X_i})$ and
\eqref{eg1-y}  also becomes   $\dot{Y}=-\nabla_Y(   \sum_i \frac{ Y_i^2}{2\Gamma_i(X)}  )$.
However, {even in this positive-definite case,} it is easy to see that
 there is no single  potential for \eqref{eg1} in the extended space $\mathcal{X}\times\mathcal{Y}$,  for whatever choice of $\Gamma$.
 This means that we do not know {\it   a priori } if  the averaging system
  is gradient or not.
The analytical form of the averaged system
\eqref{eg1-F} (with a positive-definite $D$) shows that this system  is   gradient  if and only if
the vector field $\Gamma$ is gradient, i.e., there exists a scalar function $R$ such that
\begin{equation} \label{ex1-R}
\Gamma(x) = (\Gamma_1(x),\ldots, \Gamma_n(x))=\nabla_x R(x).
\end{equation}
This suggests that the existence of an extended potential function
is a sufficient, but not  a necessary condition  for the averaged dynamics
to be gradient.

In the next, we show   the numerical results of  the MsGAD for this example.
First, we have to run  both directions $v(t)$ and $w(t)$  in our  MsGAD scheme,
because, as we  just mentioned, we do not know $DF$ is symmetric  {\it   a priori}.
To validate our result, we not only show the convergence to the saddle point,
but also compare the trajectory $x(t)$ in the MsGAD with that
of the classic  GAD applied to the known  limit equation \eqref{eg1-F}.
The parameters we used in the numerical tests are the following.
$\sigma^2=10$,
$D=\begin{bmatrix} 0.8 & -0.2  \\ -0.2 & 0.5 \end{bmatrix},
$ and for $x=(x_1,x_2)$,
$$ R(x)=\sum\limits_i \arctan(x_i-5),  ~~ \Gamma_i(x)= \left( 1+ (x_i-5)^2 \right)^{-1},~ i = 1,2.$$
The averaged equation   \eqref{eg1-F} becomes
\begin{equation}\label{eg1-ave}
\dot{\Xb} = -\nabla W(\Xb), ~~\mbox{ where } ~W(x) = x^\tr Dx- \frac{\sigma^2}{2} R(x).
\end{equation}
 $W$ has three local minima $m_1 = (0.4643, 0.6985), m_2 = (2.2038, 5.9804)$, $m_3 = (5.7109, 6.2369)$ and
  two saddle points $s_1 = (1.2841, 3.4483), s_2 = (3.5689, 6.0735), $.
 See Figure \ref{traj_Ms} below.

 In the HMM scheme of the MsGAD,  we use the forward Euler solver for the whole system.
 We take the macro-time step size $\Delta t=0.01$ for both $x^\eps$ and $v^\eps, w^\eps$
 and set $\tau=1.0$ and  $K=1$ for this example.
  The micro-time step size is ${\delta t }= \eps \times 0.01$ and the total sampling time $T=10$ is used to estimate the effective force and the Jacobi matrix.
 The initial values for $x$ are set on three local minima, respectively.
 The initial values for the directions $(v,w)$ and the fast processes are arbitrarily chosen.
 Figure \ref{traj_Ms} shows the four GAD trajectories of the $x$ component (dashed line)
 starting from three local minima.
 Depending on the initial values of $x$,  these
 four trajectories converge to the different neighboring saddle points.
 Two of them which start  from $m_2$, converge to the saddle point
 $s_1$ and $s_2$ respectively,  due to the different initial values
 for the directions $v$ and $w$.

\begin{figure}[!htb]
\begin{center}
  \includegraphics[scale=0.50]{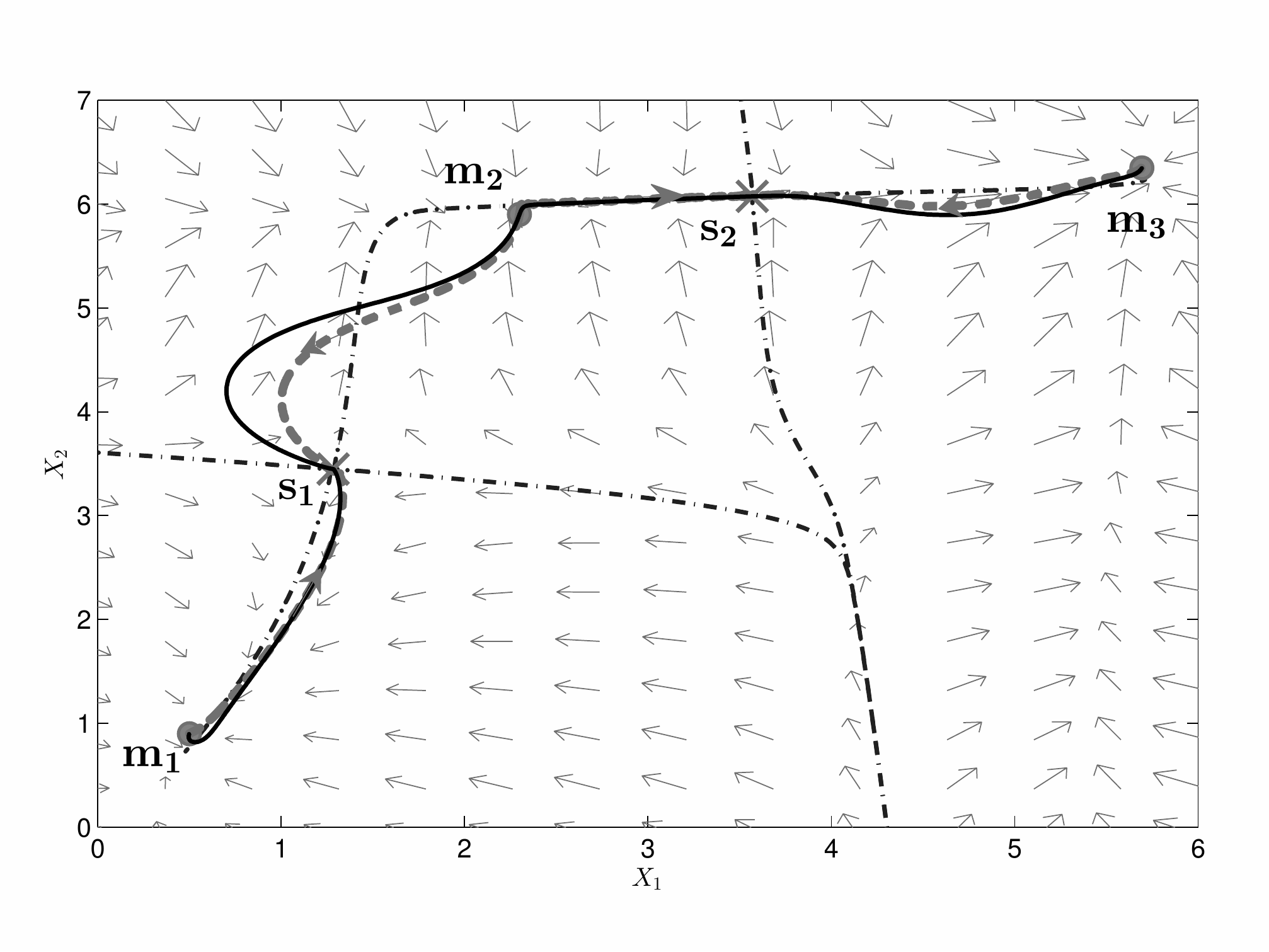}
\end{center}
 \caption{GAD trajectories from three different local minima
 ($m_1, m_2$ and $m_3$) to two different saddle points ($s_1$ and $s_2$).
The flow indicated by the  arrows corresponds to  the averaged gradient  system \eqref{eg1-ave}.
The dash-dotted curves are the stable/unstable manifolds of the two saddle points;
they determine  the basin boundaries of the three local wells.
The   thick dashed curves with arrows marked are the trajectories of MsGAD by the HMM.
As comparison, the thin solid curves are the trajectories of the GAD directly applied to
the closed form  \eqref{eg1-ave} of the averaged system.
}
\label{traj_Ms}
 \end{figure}

 {
In the end, we conduct a numerical investigation of the variance
in estimating the force $F=\bar{f}$ and the Jacobi matrix $DF=\overline{D_x f}+\Cb$.
Specifically, we look at $\var_y(f(x,y))$ and $\var_y(D_x f(x,y)+C(x,y))$,
where the variance acts   component-wisely.
The calculation shows that
$\var_y(f(x,y))=\var_y (y_i^2)=\frac{1}{2}\sigma^4 \Gamma_i(x)^2$ and
$\var_y (D_x f(x,y) + C(x,y))=\var_y(C(x,y))$.
For each component,
$\var_y(C_{ij}(x,y))=\var_y (f_i(x,y)g_j (x,y) )$,
one can find that
$f_i(x,y) g_j (x,y)= -(D_{i1}x_1 + D_{i2}x_2 + y_i^2) ( {2y_j^2}/{\sigma^2}+   \Gamma_j(x)) (x_j-5)
$
 after plugging in $\Gamma_i'(x_i) = -2(x_i-5)\Gamma^2_i(x_i)$. Thus,
 we have the analytical result that
 $\var_y(C_{ij}(x,y)) = w_{ij}(x) (x_j-5)^2$.
The expressions of  $ {w_{ij}(x)}$ are very long and there is no need to write down.
The  key observation we want to draw attention to  is that $\var_y(C_{ij}(x,y))$ is mainly dominated by the term $(x_j-5)^2$, while the variance  $\var_y(f_i(x,y))$ is proportional to $\Gamma_i(x_i)=(1+(x_i-5)^2)^{-1}$.
In Figure \ref{fig:VarDF}, we plot the variance of all components for the
force and the covariance matrix and it clearly shows that
the variance of $C_{i2} $ near $x_2\approx 5$
are distinctively different from that of  $F$ or the other two components.
Therefore, there seems no  general conclusion that the variance of the Jacobi matrix
would be larger than the force due to  the non-trivial   dependency of $\mu_x$ on $x$.
\begin{figure}[htbp]
\begin{center}
\includegraphics[width=0.65\textwidth]{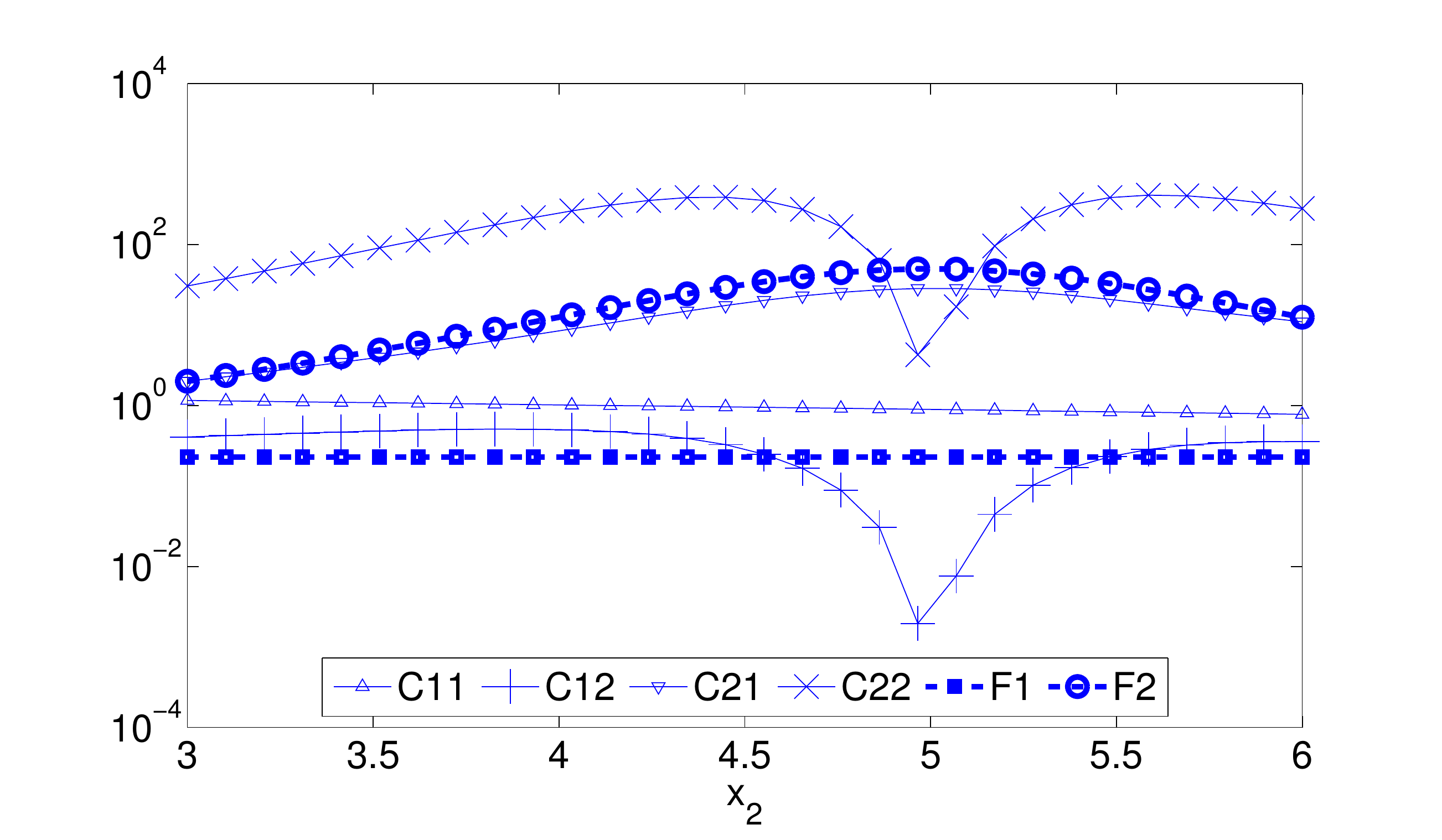}
\caption{The variances of the estimator for the effective force  $F$ and the Jacobi matrix
$DF$. In our example, the variance of $DF$ equals the variance of $C$.
The plot is the variance for each  component along a vertical line segment in $\mathcal{X}$ plane
with $x_2\in [3,6]$ and $x_1=1.2841$. The saddle point $s_1$ is on this line segment.  }
\label{fig:VarDF}
\end{center}
\end{figure}
}

\subsection{ A coupled Allen-Cahn system}
 
Our second example is  a  system of stochastic
partial differential  equations of $(u^\eps(x,t), \phi^\eps(x,t))$
in the Hilbert space $L^2([0,1])$,  satisfying      Allen-Cahn-type
equations    with   Neumann boundary condition:
\vspace{-22pt}\begin{subequations}
\begin{center}
\begin{empheq}[left=\empheqlbrace]{align}
  \partial_t u^\eps   &= \kappa^2\Delta u^\eps + u^\eps - (u^\eps)^3 + \mu\phi^\eps ,\label{eg2-u}\\
 \partial_t \phi^\eps & =   \frac{1}{\eps}[\Delta\phi^\eps - \phi^\eps + \mu u^\eps]
   + \frac{\sigma}{\sqrt{\eps}}\dot{ W}(t), \label{eg2-phi}
 \end{empheq}
\end{center}
\label{eg2}
\end{subequations}
where  $\Delta=\partial^2_x$  and
 $W(t)$ is an $L^2([0,1])$-valued Wiener process with
  {a positive-definite} 
 (spatial) covariance operator $Q$. $\dot{W}$ is white noise in time.
 $\kappa$ is the diffusion coefficient in slow dynamics
 and  $\mu$ is the coupling  constant between the slow and the fast dynamics.
$\sigma$ is the noise intensity.
 For any fixed $u^\eps = u,$ the equilibrium distribution
 of the SPDE (\ref{eg2-phi}) is the Gaussian measure
 $\mathcal{N}(\mu(I-\Delta)^{-1}u,  {\sigma^2(I-\Delta)^{-1}Q}/{2})$  on  the Cameron--Martin space. We simply choose $Q$ as identity, i.e., formally,  $\dot{W}$ is the spatio-temporal white noise. Then, it is easy to see that the averaged equation for the limit solution $\bar{u}$
 is
 \vspace{-22pt}
\begin{subequations}
\begin{center}
\begin{empheq}[left=\empheqlbrace]{align}
& \partial_t \bar{u}  = \kappa^2\Delta \bar{u} + \bar{u} - \bar{u}^3 + \mu^2 (I-\Delta)^{-1}\bar{u},\\
& \frac{\partial \bar{u}}{\partial \vec{n}}\Big|_{x=0} = 0, \quad \frac{\partial \bar{u}}{\partial \vec{n}}\Big|_{x=1} = 0.
 \end{empheq}
\end{center}
\label{eg2-F}
\end{subequations}
We can find an energy functional $U(u,\phi)$ jointly for the pair $(u,\phi)$:
\begin{equation}
U(u,\phi) = \int_\Omega \frac{\kappa^2}{2} u_x^2 + \hfour(u^2-1)^2 - \mu u \phi + \hh \phi_x^2 + \hh \phi^2 ~\dx.
\end{equation}

Since this  is a gradient system,   the corresponding MsGAD  only involves one direction
(see Section \ref{sec:grad}):
 \vspace{-22pt}\begin{subequations}
\begin{center}
\begin{empheq}[left=\empheqlbrace]{align}
  \partial_t {u}^\eps  &= -\delta_u U(u^\eps,\phi^\eps) + 2  \frac{\inpd{\delta_u U(u^\eps,\phi^\eps)}{ v^\eps}}
  { \inpd{v^\eps}{v^\eps}} v^\eps,\label{GAD1-g-u}
  \\
     \partial_t {\phi}^\eps  &=  -\frac{1}{\eps} \delta_\phi U(u^\eps,\phi^\eps)   + \frac{\sigma}{\sqrt{\eps}} \dot{W}(t),  \label{GAD1-g-phi}\\
 \partial_t{v}^\eps  & =   - \delta_u^2 U(u^\eps,\phi^\eps)v^\eps +
Cv^\eps - \alpha^\eps v^\eps, \label{GAD1-g-v}
 \end{empheq}
\end{center}
\label{GAD1-g}
\end{subequations}
 where $\delta_u U$ and  $\delta_\phi U$ are the \Fd\ of $U(u,\phi)$.
 $\delta^2_u U$ is the Hessian.
 Here $C = -  \frac{2}{\sigma^2}  {\delta_u U \otimes \delta_u U }  + \frac{2}{\sigma^2}
\overline{\delta_u U}\otimes \overline{\delta_u U}.$

\begin{figure}
 \centering
\vskip -3.5cm
 \includegraphics[scale=0.31]{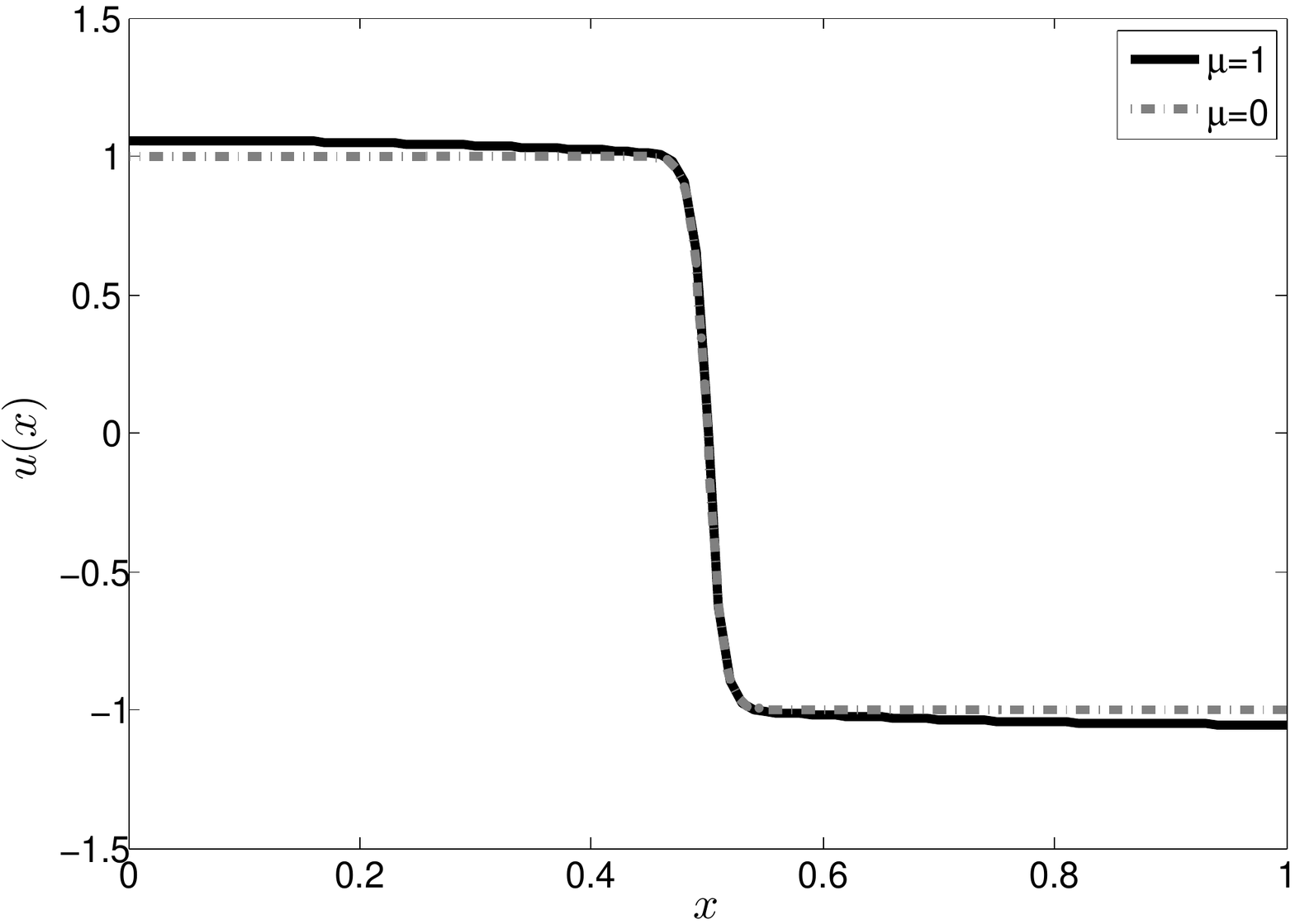}\includegraphics[scale=0.31]{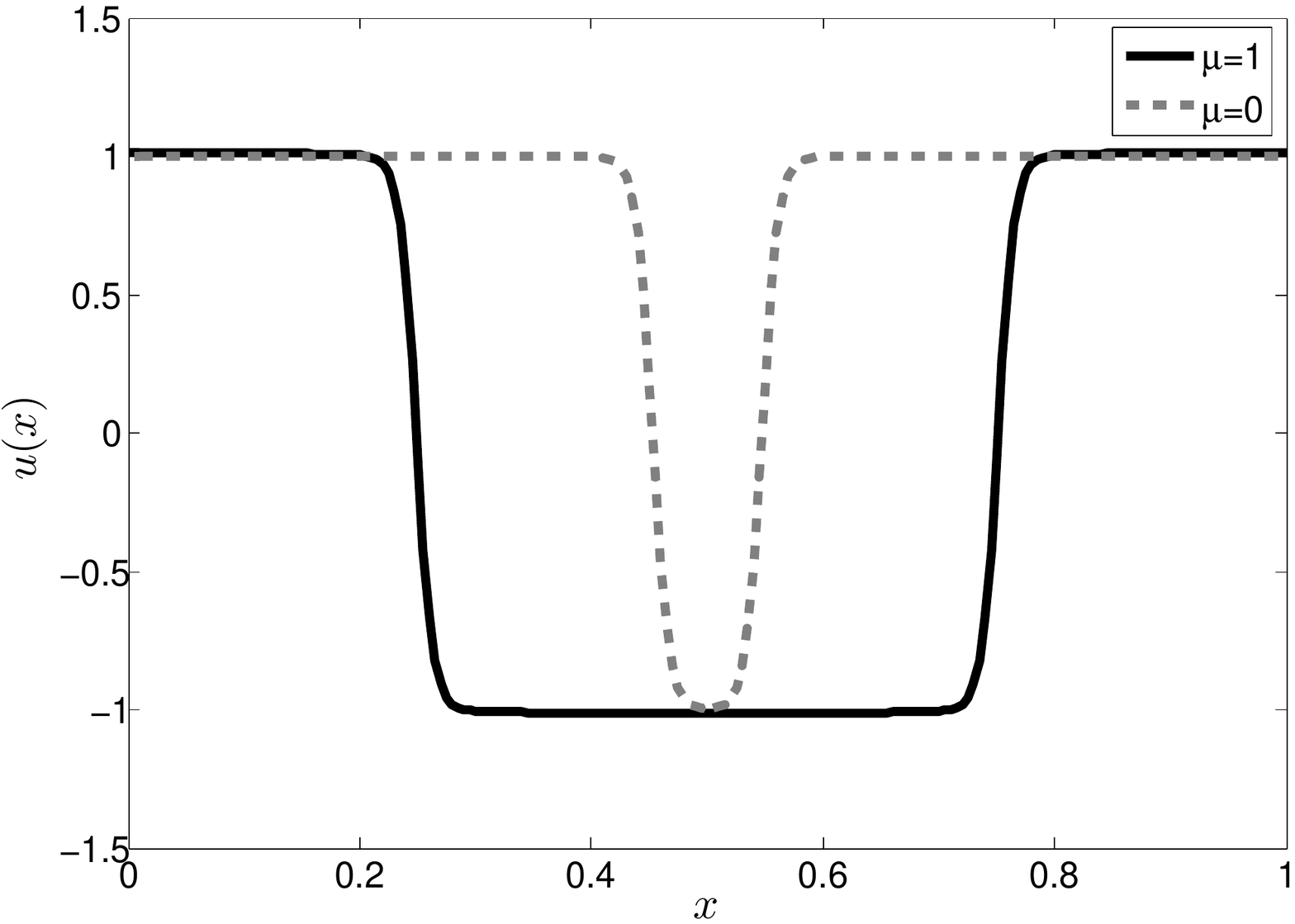}\\
  \vskip -2.4cm
   \hspace{0.2cm} (a)\hspace{6.0cm} (b)
 \vskip -0.2cm
 \caption{Two types of  index-1 saddle points for $\mu=0$ (dashed lines) and $\mu=1$
 (solid lines). The saddle point in (a) has a lower energy than the saddle point in (b).
 }
 \label{fig:ts}
 \end{figure}

In this example, the parameters are set  as $\kappa=0.01$,   $\sigma=0.3$ and 
  $\mu=1$.
 { The two local minima of the effective dynamics    are $u \equiv \pm 1$ for $\mu=0$, and $ u \equiv \pm 1.4142$ for $\mu=1$, respectively. 
The saddle points are shown in Figure \ref{fig:ts}.}

The   equation  \eqref{GAD1-g} are solved by the HMM
and the SCM.
The time-discrerization scheme for the equation \eqref{GAD1-g-u}
is a convex-splitting scheme for saddle point search
(\cite{CvxIMA2016}) to allow a   stable  time step size $\Delta t$.
The spatial discretization is the central finite difference method with the uniform mesh
grid size $\Delta x=1/200$.
The   term $ \delta_u^2 U(u,\phi) v$ in (\ref{GAD1-g}) is calculated by finite difference approximation
$$ \delta_u^2 U(u,\phi)v = \frac{1}{h} [\delta_u U(u+hv,\phi) - \delta_u U(u,\phi)]$$
with $h=0.001$.

For the seamless coupling method to solve the MsGAD \eqref{GAD1-g}, according to
Section \ref{ssec:ss}, the boosted system is
\vspace{-20pt}
 \begin{subequations}
\begin{center}
\begin{empheq}[left=\empheqlbrace]{align}
   \partial_t {u}^\eps  &= -\delta_u U(u^\eps,\phi^\eps) + 2  \frac{\inpd{\delta_u U(u^\eps,\phi^\eps)}{ v^\eps}}
  { \inpd{v^\eps}{v^\eps}} v^\eps,\label{GAD3-g-u}\\
     \partial_t {\phi}^\eps  &=  -\frac{1}{\eps\lambda} \delta_\phi U(u^\eps,\phi^\eps)   + \frac{\sigma}{\sqrt{\eps\lambda}} \dot{W}_1(t),  \label{GAD3-g-phi}\\
     \partial_t {\psi}^\eps  &=  -\frac{1}{\eps\lambda} \delta_\psi U(u^\eps,\psi^\eps)   + \frac{\sigma}{\sqrt{\eps\lambda}} \dot{W}_2(t),  \label{GAD3-g-psi}\\
  \partial_t{v}^\eps &= - \delta_u^2 U(u^\eps,\phi^\eps)v^\eps +
 f(u^\eps,\phi^\eps) \inpd{g(u^\eps,\phi^\eps)} {v^\eps}  \notag
 \\
 &~~\quad
 -f(u^\eps,\psi^\eps) \inpd{g(u^\eps,\phi^\eps)} {v^\eps}
  - \alpha^\eps v^\eps, \label{GAD3-v}
 \end{empheq}
\end{center}
\label{GAD3}
\end{subequations}
where $W_1$ and $W_2$ are two iid copies of $W$.
$\inpd{\cdot}{\cdot}$ is the inner product in $L^2([0,1])$.

Define the errors as follows,
\[ \mbox{err}_H(t) = \|u_H(x,t)-u^*(x)\|_{L^2},~~~\quad
\mbox{err}_S (t) = \|u_S(x,t)-u^*(x)\|_{L^2},
\]
where $u_H(x,t)$ represents the result of \eqref{GAD1-g} solved  by
the HMM and $u_S(x,t)$ is the result of \eqref{GAD3} solved  by the SCM.
 The
true solution  $u^*(x)$ is  obtained from the classic GAD applied to  the closed form of the averaged system
\eqref{eg2-F} with a very fine time step size and a sufficiency small tolerance.

To test our HMM and SCM, 
we calculate the saddle point with the lowest energy, i.e., the profile in the left panel of Figure \ref{fig:ts}.
The initial guess is   the function $\cos(\pi x)$.

 \begin{figure}[htbp]
 \centering
 \includegraphics[width=0.53\textwidth]{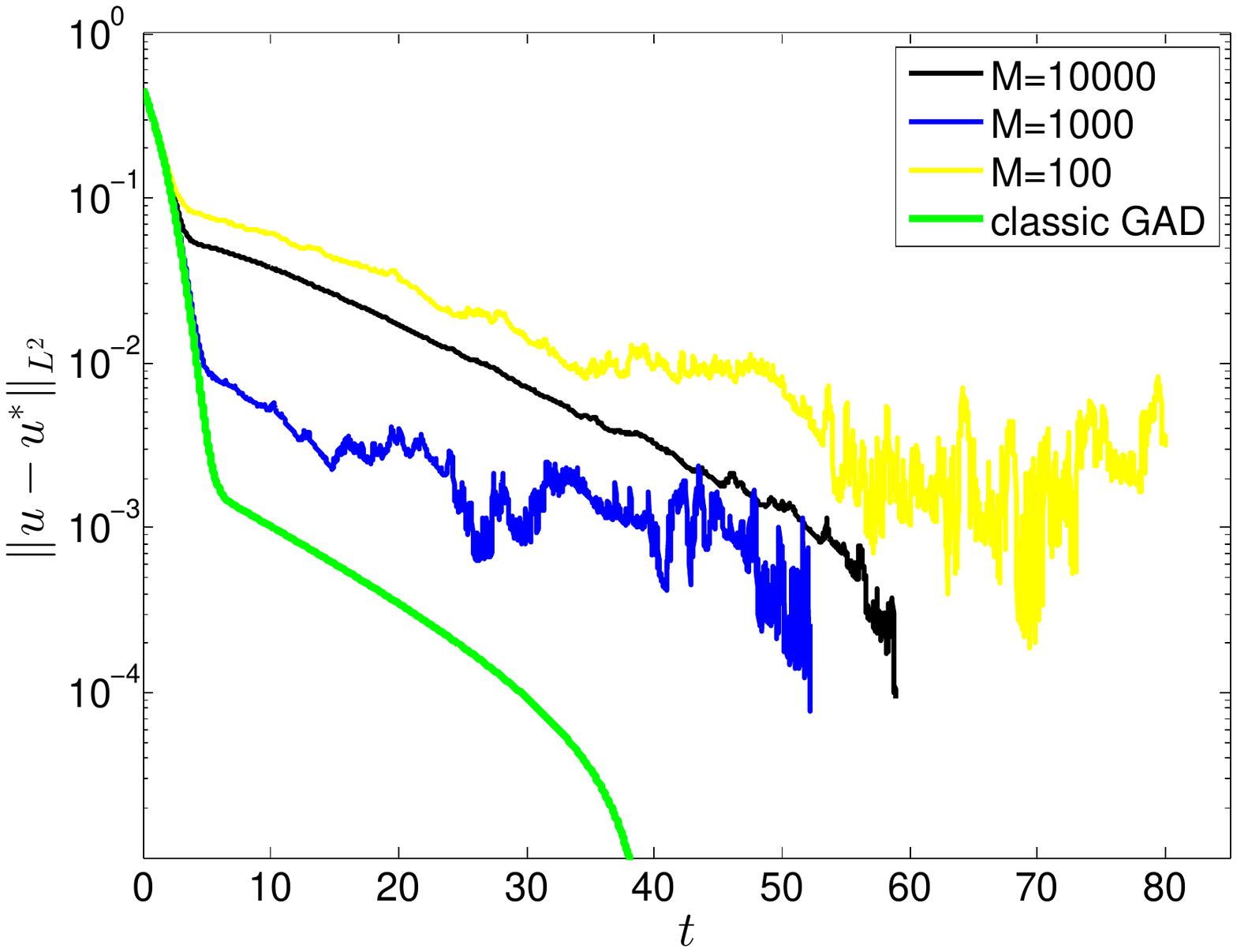}
 \hspace{-10mm}\includegraphics[width=0.53\textwidth]{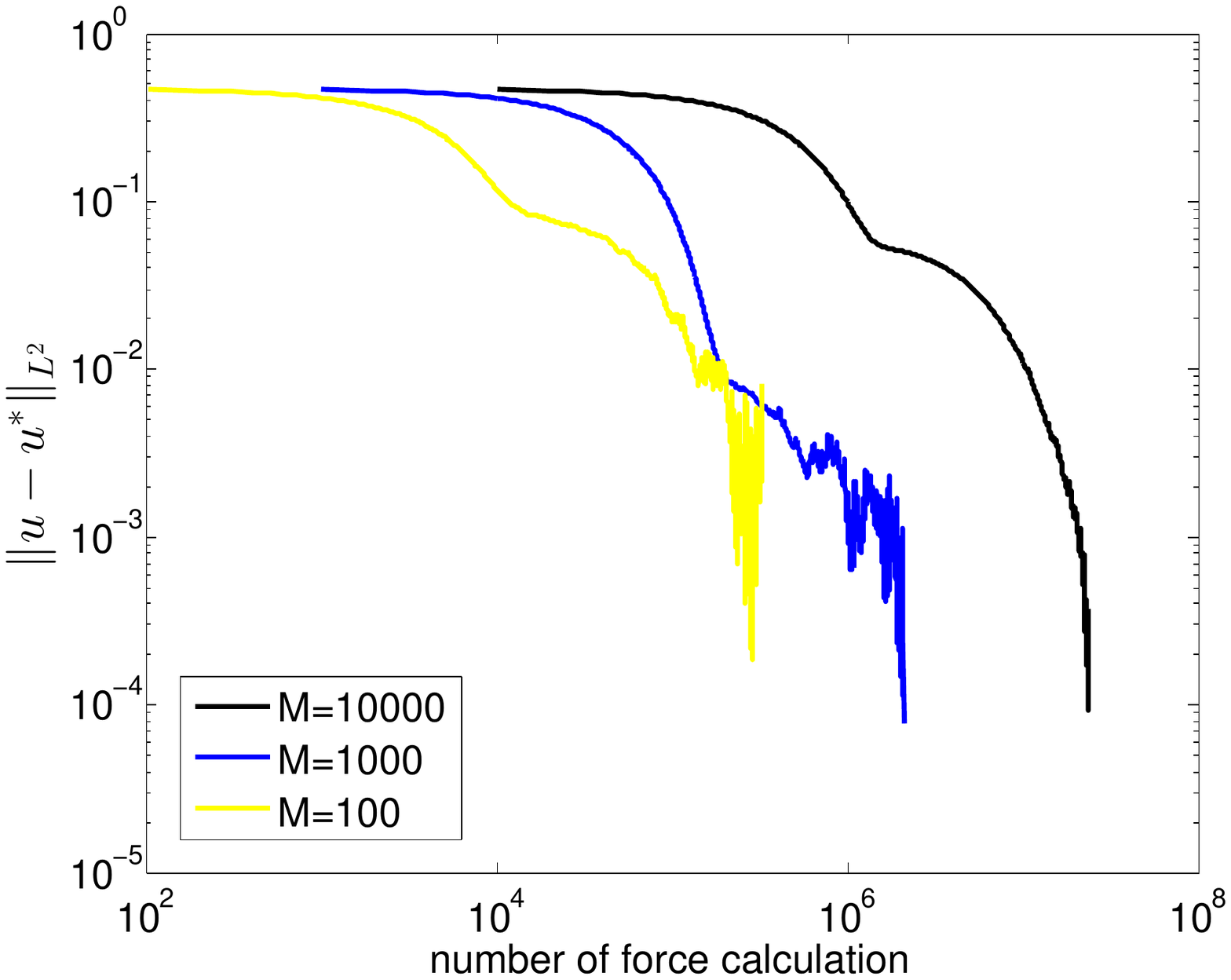}\\
  \vskip -2.6cm
    (a)\hspace{5.0cm} (b)
 \vskip -0.2cm
\caption{(a) shows the evolution of the error with the physical time for  the classic GAD applied to \eqref{eg2-F} and for the HMM with various sample sizes $M=1000$ and $10000$. (b) is the  decay of the errors with the number of force calculations.}\label{HMM_error}
 \end{figure}

In the HMM, the macro-time step size for $u^\eps$ and $v^\eps$
is $\Delta t = 0.025$ and  the micro-time step size for $\phi^\eps$
is $\delta t=0.01\times \eps$.
To demonstrate  the effect of the sample size $M$ in the HMM, we use different  $M$ and run the simulation
up to the error $10^{-4}$ and plot in  Figure \ref{HMM_error}  the errors w.r.t.
 the time  as well as w.r.t.
   the number of force calculations (including the cost in calculating the direction variables). 
It is observed that  a smaller sample size leads to a relatively larger fluctuation and the fluctuation appears earlier in time.  The reduced cost from a smaller $M$ is quite obvious in  the subfigure (b). 
One can  expect a further cost reduction if   $M$ decreases further;  but to attain the error as small as $10^{-4}$,
one should either switch to a larger $M$ when approaching the saddle point
or average  the  output in time in order to reduce the variance.

Next, we present the details of the  SCM and its adaptive version.
In the standard SCM, the constant time step size $\wt{\Delta }t = 1.0\times 10^{-2}$ is applied  for all components and $\eps = 1.0\times 10^{-4}, \lambda =10$ (or effectively, $\eps^{\prime}=\eps\lambda=1.0\times 10^{-3}$).
In the adaptive SCM,  the same values of   $\wt{\Delta }t$ and $\eps'$ are used up to the time $t_a$, where  the error $\mbox{err}_S$ attains  a given threshold $10^{-2}$.
From the time $t_a$,
the adaptive time step size $\wt{\Delta} t$ and the adaptive $\eps^{\prime}$ are then applied. 
Specifically,  
we use $\wt{\Delta} t_k = k^{-1/2} \wt{\Delta} t$ and $\eps'_k = k^{-1/2} \eps^{\prime}$
where the subindex $k$ counts the number of time steps starting  from $t_a$.
As stated in Section \ref{Sec2_3},
the final output of the SCM and adaptive SCM are the time-averaged solution
 $ \tilde{u}(x,t) \doteq \frac{1}{t-t_0}\int_{t_0}^{t} u(x,s)\d s$,
 where $t_0=10$ in this example.
This continuous-time integral can be easily implemented 
 on the discrete time points in an iterative way.

  \begin{figure}[htbp]
 \centering
\includegraphics[width=0.51\textwidth]{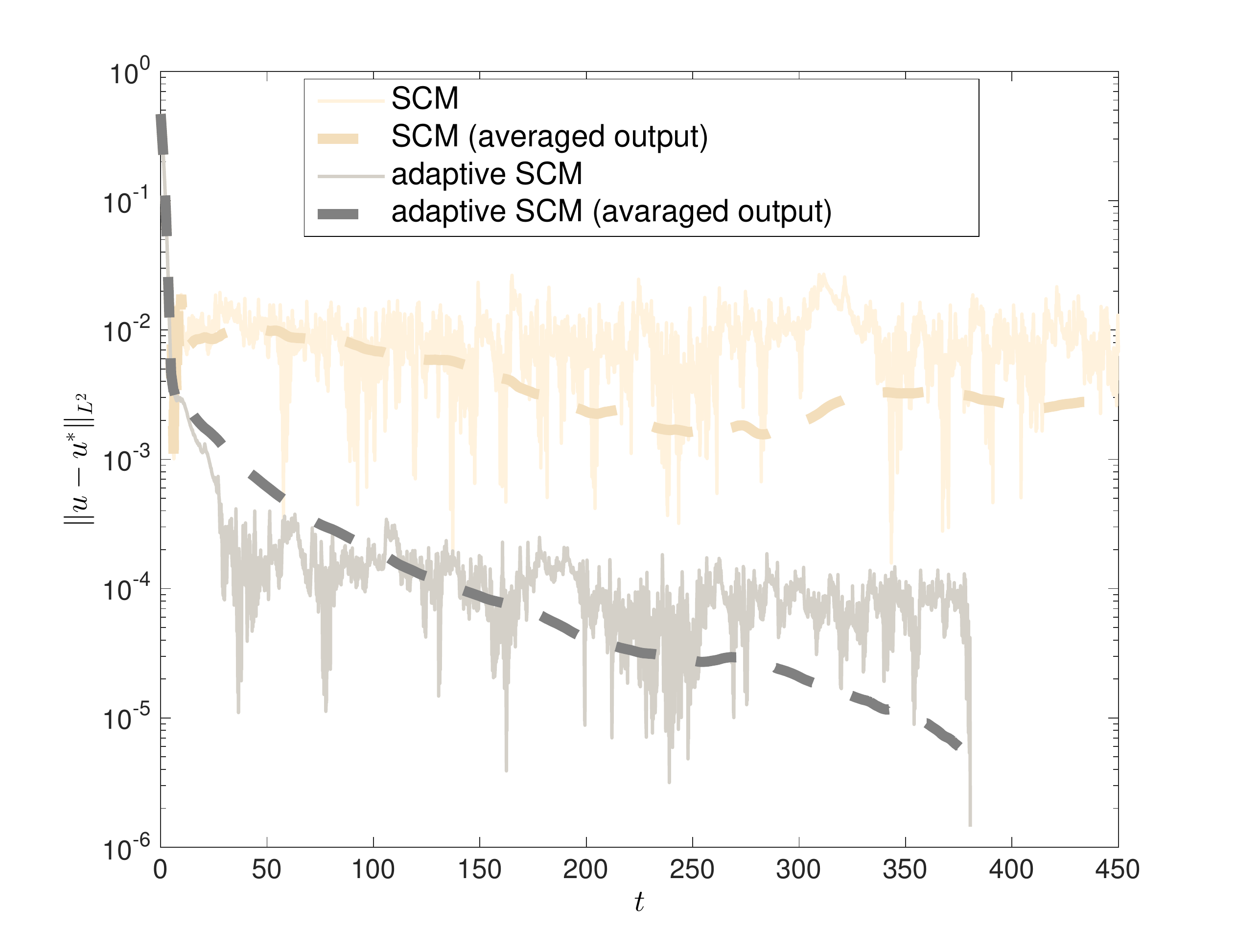}
\hspace{-15pt}
\includegraphics[width=0.51\textwidth]{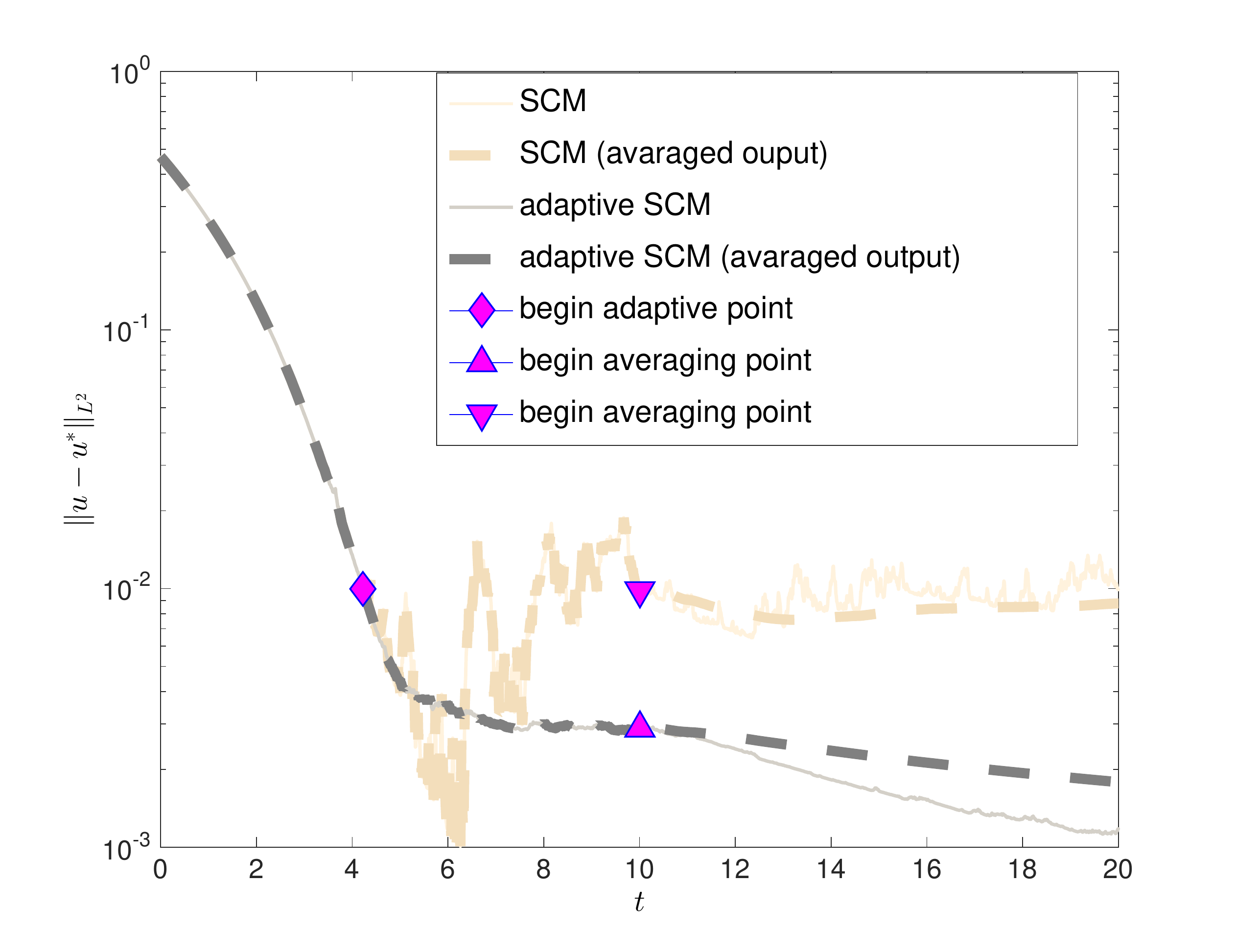}
\caption{ The evolution of the error with the time in the seamless coupling method.  The right panel is the same plot in the early stage before time $t=20$.
The solutions are  averaged w.r.t. time after the burn-in time $t_0=10$, which are shown by the two triangular symbols in the plot.
 The resulted new outputs are shown in dashed thick lines.
The  value of $\eps'$ and the step size in the adaptive SCM start to  decay when the error decreases to $0.01$,  which is shown by the diamond symbol in the right plot.
The ``SCM'' means the standard SCM with a fixed $\eps'$ and a constant time step size.  }\label{scm_error}
 \end{figure}

\begin{figure}[htbp]
 \centering
\includegraphics[scale=0.35]{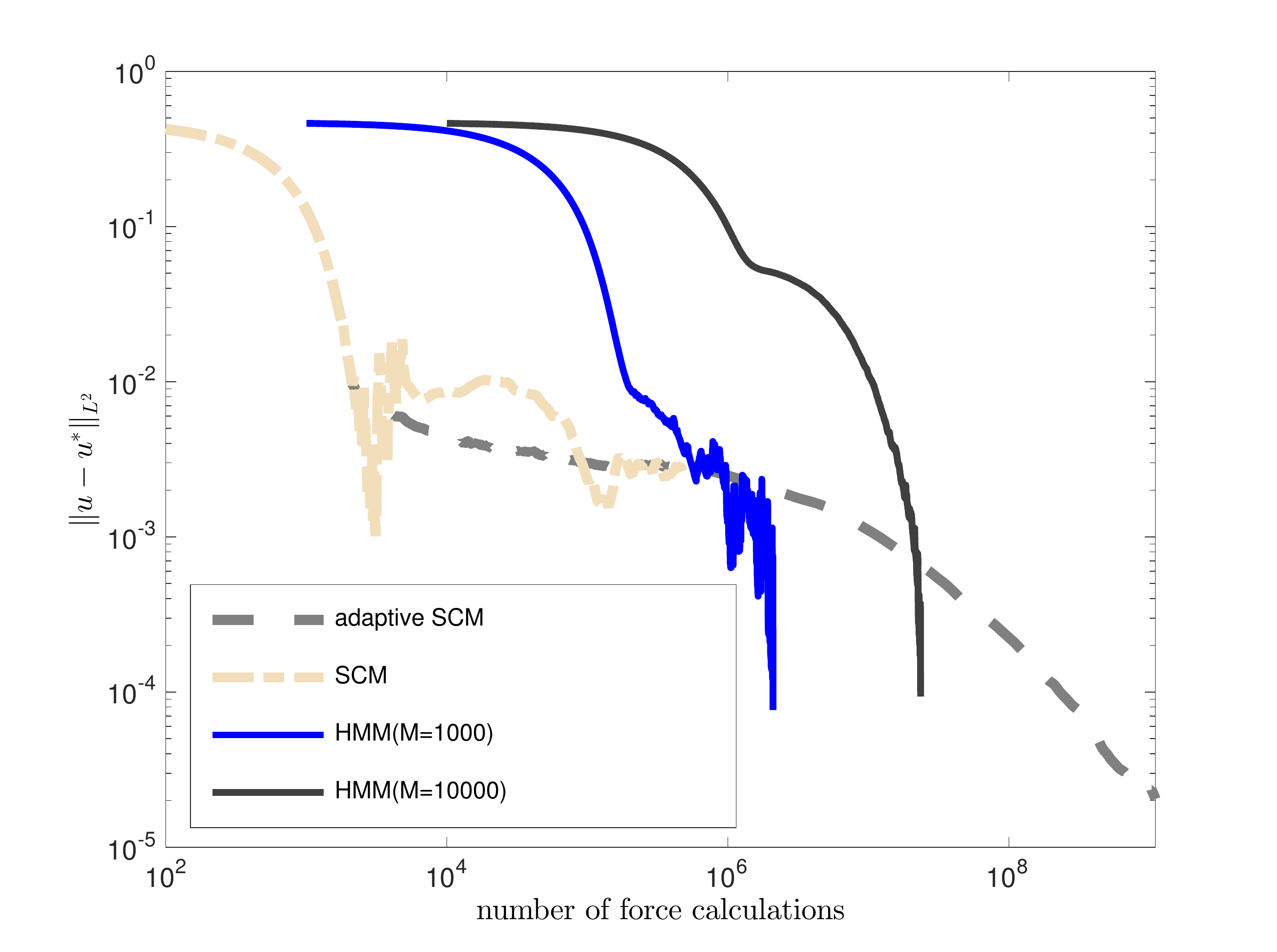}
\caption{The illustration of the computational costs for the SCM, the adaptive version of SCM (both with averaged output) and the HMM.
The plot is in log-log scale.}\label{fig:force}
 \end{figure}

 Figure \ref{scm_error} shows the effectiveness of the strategy of  decreasing  $\eps'$ at time marches. For a fixed $\eps'=10^{-3}$, the error of the SCM can not reach the level of $10^{-3}$.
 This discrepancy quickly shrinks when we start to decrease the value of $\eps'$ gradually. 
 The error eventually can decrease to the level $10^{-5}$.
 However, this comes with a price: as $\eps'$ shrinks, the problem becomes stiffer 
      and the step size $\wt{\Delta} t$ has to decrease in the same speed.
   The log-log plot in  Figure \ref{fig:force} shows the fast growth of the computational cost
   when $\eps'$ decreases in the adaptive SCM.  Thus,  the SCM is not
   desired for a tiny $\eps'$, though it can quickly settle down to a region of importance with cheap cost.
   Figure  \ref{fig:force} also illustrates the performance of the HMM and the SCM.
   Note that since the time step sizes used here  in HMM and  adaptive  SCM are different, this figure 
   is only a qualitative picture. But what we can learn from the above numerical explorations
   for this example  is that the best strategy in practice perhaps is 
   to start with a boosted value  $\eps'$ in the adaptive SCM
and then switch to the HMM with an adaptive choice of $M$.
The optimal choices of the underlying parameters are of both practical importance and theoretic interests
and are left for future study.

\section{Conclusion}
\label{sec:end}
We have proposed the multiscale method based on the GAD 
 in  slow-fast systems.
In order to calculate the saddle point of effective dynamics,
 we derived a new slow-fast system, the MsGAD, whose effective
dynamics is consistent with the GAD for the  effective dynamics.
By applying the multiscale numerical method such as the HMM or the
seamless coupling method,  we   efficiently compute   saddle points
on   averaged dynamics.
{Some adaptive techniques are presented and illustrated by the numerical experiment. 
}

 \bigskip
 \section*{Appendix: Proof of \eqref{B}.}
To show that $B$ is  the Jacobi matrix $D F$, we just need to prove that
\begin{equation*}
\begin{split}
\int (f(x,y)-F(x)) \otimes  \nabla_x \log \rho(x,y) \rho(x,y)(\dy)
\\=   \int  \left(  \int_0^\infty \nabla_y \e^y f(x, \wt{Y}^x_\tau)\d \tau
\right) D_x b(x,y)  \rho(x,y)\dy.
\end{split}
\end{equation*}
To show this matrix equality,  we just need to show
\[
\begin{split}
&\int (f(x,y)-F(x))  \inpd{  \nabla_x   \rho(x,y) }{e} \dy
\\= &  \int  \left(  \int_0^\infty \nabla_y \e^y f(x, \wt{Y}^x_\tau)\d \tau
\right) (D_x b(x,y) )e\, \rho(x,y) \dy
\end{split}
\]
for   arbitrary vector $e$.
Write  the infinitesimal perturbation
$\tm :=
  \rho(x+he,y)-\rho(x,y)   \approx h\nabla_x \rho(x,y) \cdot e $
 and   $\wt{b}(x,y)=b(x+he,y)-b(x,y)\approx h ( D_x b(x,y) )e$ for $h\ll 1$,
and assume $e$ is independent of $x$, then
  the   lemma below tells us that
\[ \int (f(x,y)-F(x))  \inpd{  \nabla_x   \rho(x,y) }{e} \dy=   -\int
  \inpd{\nabla_y u} {(D_x b(x,y)) e}  \rho(x,y) \dy,
\]
where $u$ satisfies
$\Lo_y u(x,y) = f(x,y) - F(x) $.
By the Feymann-Kac formula, we have the representation
$u(x,y) =  \int_0^\infty\e^y [-f(x,\wt{Y}^x_t)+F(x)]\dt $,
so $\nabla_y u = - \int_0^\infty \nabla_y \e^y [f(x,\wt{Y}^x_t)]\dt$.
This completes our proof.

\begin{lemma}
If  $\rho(y)$ is the unique equilibrium probability density function  of the SDE
\[  \d Y = b(Y) \dt + \sigma(Y) \d W,\]
 $\Lo$ is the infinitesimal  generator, and the infinitesimal perturbation is applied for the drift term:
$b \to b+\wt{b}$ and the diffusion term $a\doteq \sigma\sigma^\tr \to a + \wt{a}$,
where $\wt{b}$ and $\wt{a}$ are small terms,
and let   $\Theta\doteq \int \theta(y) \rho(y)\dy$, then
the perturbation of $\Theta$ is
\[ \delta \Theta \doteq \int \theta(y) \tm(y)\dy =- \int  \wt{b} \cdot ( \rho \nabla u) (y) \dy  - \frac12 \int \wt{a} : \rho \nabla \nabla u (y)\dy,\]
where $u$ is the solution of the adjoint equation $\Lo u = \theta(y) - \Theta$ and decays to $0$
at infinity.
\end{lemma}
{\it Proof:}
The infinitesimal generator is $\Lo = b(y) \cdot \nabla + \frac12  a(y) : \nabla \nabla$.
The density $\rho$ satisfies the equation $\Lo^* \rho=0$, where
$\Lo^*$ is the adjoint of $\Lo$. By linearizing the perturbed equation,
we have $\Lo^*\wt{\rho} = \nabla \cdot (\rho \wt{b})  - \nabla \nabla : (\rho \wt{a})$.
Multiply this equation by $u$, then from the integration by parts, we have
$\delta \Theta = \int \wt{\rho}(y)\theta(y) \dy = \int \wt{\rho}(y)\Lo u  \dy
= - \int  \wt{b} \cdot ( \rho \nabla u) (y) \dy  - \frac12 \int \wt{a} : \rho \nabla \nabla u (y)\dy.$

\bibliography{./MsGAD}

\bibliographystyle{siam} 

\end{document}